\documentclass[10pt]{article}

\usepackage{amsmath,amsxtra,amssymb,latexsym, amscd,amsthm}
\usepackage{here}
\usepackage{graphpap}
\usepackage{graphicx}
\usepackage{color}
\usepackage{ marvosym }

\newtheorem{theorem}{Theorem}[section]

\newtheorem{lemma}[theorem]{Lemma}
\newtheorem{proposition}[theorem]{Proposition}
\theoremstyle{definition}
\newtheorem{definition}[theorem]{Definition}
\newtheorem{remark}[theorem]{Remark}

\newtheorem{example}[theorem]{Example}

\numberwithin{equation}{section}

\def\divv{\text{div}}
\def\<{\langle}
\def\>{\rangle}
\def\R{\mathbb{R}}
\def\N{\mathbb{N}}
\def\limn{\lim_{n\to\infty}}
\def\limsupn{\limsup_{n\to\infty}}
\def\liminfn{\liminf_{n\to\infty}}
\definecolor{purple}{rgb}{0.4, 0.0, 0.4}

\begin{document}

\begin{center}
{\bf \Large Identifying conductivity in electrical impedance tomography \\with total variation regularization}
\end{center}

\vspace{0.5cm}

\centerline {\bf Michael Hinze$^a$\let\thefootnote\relax\footnote{Email: michael.hinze@uni-hamburg.de,~ Barbara.Kaltenbacher@aau.at,~ quyen.tran@uni-hamburg.de} $\mathbf{\cdot}$ Barbara Kaltenbacher$^b$\let\thefootnote\relax\footnote{M. Hinze gratefully acknowledges support of the Lothar Collatz Center for Computing in Science at the University of Hamburg}\let\thefootnote\relax\footnote{B. Kaltenbacher gratefully acknowledges support of the Austrian Wissenschaftsfonds through grant FWF I2271 entitled "Solving inverse problems without forward operators"}\let\thefootnote\relax\footnote{T.N.T. Quyen gratefully acknowledges support of the Alexander von Humboldt-Foundation} $\mathbf{\cdot}$ Tran Nhan Tam Quyen$^{a,c}$\footnote{$^c$Author to whom any correspondence should be addressed}}

{$^a$University of Hamburg, Bundesstra\ss{}e 55, D-20146 Hamburg, Germany\\
$^b$Alpen-Adria-Universit\"at Klagenfurt, Universit\"atsstra\ss{}e 65-67, A-9020 Klagenfurt, Austria

%\vspace{0.5cm}
\centerline{\rule{8cm}{.1pt}}

{\small {\bf Abstract:} In this paper we investigate the problem of identifying the conductivity in electrical impedance tomography
from one boundary measurement. A variational method with total variation regularization is here proposed to tackle this problem. We discretize the PDE as well as the conductivity with piecewise linear, continuous finite elements. We prove the stability and convergence of this technique. For the numerical solution we propose a projected Armijo algorithm. Finally, a numerical experiment is presented to illustrate our theoretical results.}

{\small {\bf Key words and phrases:} Conductivity identification, electrical impedance tomography, total variation regularization, finite element method, Neumann problem, Dirichlet problem, ill-posed problems.}

{\small {\bf AMS Subject Classifications:} 65N21; 65N12; 35J25; 35R30}

\section{Introduction}

Let $\Omega$ be an open bounded connected domain in $ \R^d, d\in\{2,3\}$ with polygonal boundary $\partial \Omega$ and $f\in {H^1(\Omega)}^* := H^{-1}(\Omega)$ be given. We consider the following elliptic boundary value problem
\begin{align}
-\nabla \cdot \big(q \nabla \Phi \big) &= f \mbox{~in~} \Omega,  \label{17-5-16ct1}\\
q \nabla \Phi \cdot \vec{n} &= j^\dag \mbox{~on~} \partial\Omega \mbox{~and~} \label{17-5-16ct2}\\
\Phi &= g^\dag \mbox{~on~} \partial\Omega, \label{17-5-16ct3}
\end{align}
where $\vec{n}$ is the unit outward normal on $\partial\Omega$.

The system \eqref{17-5-16ct1}--\eqref{17-5-16ct3} is overdetermined, i.e. if the Neumann and Dirichlet boundary conditions $j^\dag \in H^{-1/2}(\partial\Omega) := {H^{1/2}(\partial\Omega)}^*, ~g^\dag \in H^{1/2}(\partial\Omega)$ and the {\it conductivity}
\begin{align}\label{22-3-15ct3}
q\in \mathcal{Q} := \left\{ q \in L^{\infty}(\Omega)
~\big|~ \underline{q} \le q(x) \le \overline{q} \mbox{~a.e. in~} \Omega\right\}
\end{align}
are given, then there may be no $\Phi$ satisfying this system. Here $\underline{q}$ and $\overline{q}$ are some given positive
constants.

In this paper we assume that the system is consistent and our aim is to identify the conductivity $q^\dag\in \mathcal{Q}$ and the electric potential $\Phi^\dag \in H^1(\Omega)$ in the system \eqref{17-5-16ct1}--\eqref{17-5-16ct3} from current and voltage i.e., Neumann and Dirichlet measurements at the boundary $\left(j_\delta,g_\delta \right)\in H^{-1/2}(\partial\Omega) \times H^{1/2}(\partial\Omega) $ of the exact $\big(j^\dag,g^\dag\big)$ satisfying
$$\big\|j_\delta-j^\dag\big\|_{H^{-1/2}(\partial\Omega)} + \big\|g_\delta-g^\dag\big\|_{H^{1/2}(\partial\Omega)} \le \delta \mbox{~with~} \delta>0.$$
Note that using the $H^{-1/2}(\partial\Omega) \times H^{1/2}(\partial\Omega)$ topology for the data is natural from the point of view of solution theory for elliptic PDEs but unrealistic with regard to practical measurements. We will comment in this issue in Remark \ref{rem:boundarydata} below.

For the purpose of conductivity identification --- a problem which is very well known in  literature and practice as electrical impedance tomography EIT, see below for some references --- we simultaneously consider the Neumann problem
\begin{align}
-\nabla \cdot (q\nabla u) = f \mbox{~in~} \Omega \mbox{~and~} q\nabla u\cdot \vec{n} = j_\delta \mbox{~on~} \partial\Omega \label{9-6-16ct2}
\end{align}
and the Dirichlet problem
\begin{align}
-\nabla \cdot (q\nabla v) = f \mbox{~in~} \Omega \mbox{~and~} v = g_\delta \mbox{~on~} \partial\Omega \label{9-6-16ct3}
\end{align}
and respectively denote by $\mathcal{N}_q j_\delta$, $\mathcal{D}_q g_\delta$ the unique weak solutions of the problems \eqref{9-6-16ct2}, \eqref{9-6-16ct3}, which depend nonlinearly on $q$, where $\mathcal{N}_q j_\delta$ is normalized with vanishing mean on the boundary.
We adopt the variational approach of Kohn and Vogelius in \cite{Kohn_Vogelius1,Kohn_Vogelius11,Kohn_Vogelius2} to the identification problem. In fact,
for estimating the conductivity $q$ from
the observation $\left(j_\delta,g_\delta \right)$ of the exact data $\big(j^\dag,g^\dag\big)$, we use the functional
$$
\mathcal{J}_\delta(q) := \int_\Omega q\nabla\left(\mathcal{N}_q j_\delta - \mathcal{D}_q g_\delta\right) \cdot \nabla\left(\mathcal{N}_q j_\delta - \mathcal{D}_q g_\delta\right)dx.
$$
For simplicity of exposition we restrict ourselves to the case of just one Neumann-Dirichlet pair, while the approach described here can be easily extended to multiple measurements $\left(j_\delta^i,g_\delta^i \right)_{i=1,\ldots,I}$, see also  Example \ref{ex_mult} below. It is well-known that such a finite number of boundary data in general only allows to identify  conductivities taking finitely many different  values in the domain $\Omega$, see, e.g., \cite{Alessandrini95}.

Indeed, we are interested in estimating such piecewise constant conductivities and therefore use total variation regularization, i.e., we consider the minimization problem
\begin{align}\label{11-7-16ct1}
\min_{q\in\mathcal{Q}_{ad}} \mathcal{J}_\delta(q) + \rho \int_\Omega \left| \nabla q\right|,
\end{align}
where
$\mathcal{Q}_{ad} := \mathcal{Q} \cap BV(\Omega)$
is the admissible set of the sought conductivities, $BV(\Omega)$ is the space of all functions with bounded total variation (see \S \ref{notation} for its definition) and $\rho >0$ is the regularization parameter, and consider a minimizer $q_{\rho,\delta}$ of \eqref{11-7-16ct1} as reconstruction.

For each $q\in\mathcal{Q}$ let $\mathcal{N}^h_qj_\delta$ and $\mathcal{D}^h_q g_\delta$ be corresponding approximations of $\mathcal{N}_qj_\delta$ and $\mathcal{D}_qg_\delta$ in the finite dimensional space $\mathcal{V}^h_1$ of piecewise linear, continuous finite elements and $q^h_{\rho,\delta} $ denote a minimizer of the discrete regularized problem corresponding to \eqref{11-7-16ct1}, i.e. of the following minimization problem
\begin{align}\label{16-6-16ct3}
\min_{q \in \mathcal{Q}^h_{ad}} \int_\Omega q\nabla\left(\mathcal{N}^h_q j_\delta - \mathcal{D}^h_q g_\delta\right) \cdot \nabla\left(\mathcal{N}^h_q j_\delta - \mathcal{D}^h_q g_\delta\right)dx + \rho \int_\Omega \sqrt{\left|\nabla q\right|^2+\epsilon^h}
\end{align}
with $\mathcal{Q}^h_{ad} := \mathcal{Q}_{ad} \cap \mathcal{V}^h_1$ and $\epsilon^h$ being a positive functional of the mesh size $h$ satisfying $\lim_{h\to0} \epsilon^h =0$.

In Section \ref{Stability} we will show  the stability of approximations for fixed positive $\rho$. Furthermore as $h,\delta \to 0$ and with an appropriate a priori regularization parameter choice $\rho=\rho(h,\delta)$, there exists a subsequence of $\big(q^h_{\rho,\delta}\big)$ converging in the $L^1(\Omega)$-norm to a total variation-minimizing solution $q^\dag$ defined by
$$q^\dag \in \arg \min_{\big\{q\in \mathcal{Q}_{ad} ~|~ \mathcal{N}_qj ^\dag = \mathcal{D}_q g^\dag \big\}} \int_\Omega | \nabla q|.$$
In particular, if $q^\dag$ is uniquely defined, then this convergence holds for the whole sequence $\big(q^h_{\rho,\delta}\big)$.
The corresponding state sequences $\Big(\mathcal{N}^h_{q^h_{\rho,\delta}} j_\delta\Big)$ and $\Big(\mathcal{D}^h_{q^h_{\rho,\delta}} g_\delta\Big)$ converge in the $H^1(\Omega)$-norm to $\Phi^\dag = \Phi^\dag(q^\dag,j^\dag,g^\dag)$ solving the system \eqref{17-5-16ct1}-\eqref{17-5-16ct3}. Finally, for the numerical
solution of the discrete regularized problem \eqref{16-6-16ct3}, in Section \ref{11-7-16ct2} we employ a projected Armijo algorithm. Numerical results show the  efficiency of the proposed method and illustrate our theoretical findings.

\medskip

We conclude this introduction with a selection of references from the vast literature on EIT, which has evolved to a highly relevant imaging and diagnostics tool in industrial and medical applications and has attracted great attention of many scientists in the last few decades.

To this end, for any fixed $q\in\mathcal{Q}$ we define the {\it Neumann-to-Dirichlet map} $\Lambda_q: H^{-1/2}(\partial\Omega)
\rightarrow H^{1/2}(\partial\Omega)$, by
\begin{align*}
j\mapsto \Lambda_qj := \gamma \mathcal{N}_qj,
\end{align*}
where $\gamma : H^1(\Omega) \to H^{1/2}(\partial\Omega)$
is the Dirichlet trace operator.
Calder\'on in 1980 posed the question whether an unknown conductivity distribution inside a domain can be determined from an infinite number of boundary observations, i.e. from the Neumann-to-Dirichlet map $\Lambda_q$:
\begin{align}\label{12-7-16ct1}
p, q \in \mathcal{Q} \subset L^\infty(\Omega) \mbox{~with~} \Lambda_p = \Lambda_q
\quad \Rightarrow \quad p=q \ ?
\end{align}
Calder\'on did not answer his question \eqref{12-7-16ct1}; however, in \cite{calderon} he proved that the problem linearized at constant conductivities has a unique solution. In dimensions three and higher Sylvester and Uhlmann \cite{sylvester} proved the unique identifiability of a $C^\infty$-smooth conductivity. P\"aiv\"arinta {\it el al.} \cite{paivarinta} and Brown and Torres \cite{brown} established uniqueness in the inverse conductivity problem for $W^{3/2,p}$-smooth conductivities with $p=\infty$ and $p>2d$, respectively. In the two dimensional setting, Nachman \cite{nachman} and Brown and Uhlmann \cite{brown2} proved uniqueness results for conductivities which are in $W^{2,p}$ with $p>1$ and $W^{1,p}$ with $p>2$, respectively. Finally, in 2006 the question \eqref{12-7-16ct1} has been answered to be positive by Astala and P\"aiv\"arinta \cite{astala} in dimension two. For surveys on the subject, we refer the reader to \cite{borcea,cheney,dijkstra,mueller,uhlmann} and the references therein.

Although there exists a large number of papers on the numerical solution of the inverse problems of EIT, among these also papers considering the Kohn-Vogelius functional (see, e.g., \cite{Knowles1998,KohnMcKenney90}) and total variation regularization (see, e.g., \cite{Dobson96,Niinimaekietal16}), we have not yet found investigations on the discretization error in a combination of both functionals for the fully nonlinear setting, a fact which motivated the research presented in this paper.

Throughout the paper we use the
standard notion of Sobolev spaces $H^1(\Omega)$, $H^1_0(\Omega)$, $W^{k,p}(\Omega)$, etc from, for example, \cite{adams}. If not stated otherwise we write
$\int_\Omega \cdots$ instead of $\int_\Omega \cdots dx$.

\section{Problem setting and preliminaries}

\subsection{Notations}\label{notation}

Let us denote by
$$\gamma : H^1(\Omega) \to H^{1/2}(\partial\Omega)$$
the continuous Dirichlet trace operator while
$$\gamma^{-1} : H^{1/2}(\partial\Omega) \to H^1(\Omega)$$
is the continuous right inverse operator of $\gamma$, i.e. $ (\gamma\circ \gamma^{-1}) g =g$ for all $g\in H^{1/2}(\partial\Omega)$.
With $f\in H^{-1}(\Omega)$ (with a slight abuse of notation) in \eqref{17-5-16ct1} being given, let us denote
$$c_f := (f,1),$$
where the expression $(f,\varphi)$ denotes the value of the functional $f\in H^{-1}(\Omega)$ at $\varphi\in H^1(\Omega)$. We also denote
$$H^{-1/2}_{-c_f}(\partial\Omega) := \left\{ j \in H^{-1/2}(\partial\Omega) ~\big|~ \langle j,1 \rangle = -c_f\right\},$$
where the notation $\left\langle j,g\right\rangle$ stands for the value of the functional $j\in H^{-1/2}(\partial\Omega)$ at $g\in H^{1/2}(\partial\Omega)$. Similarly, we denote
$$H^{1/2}_\diamond(\partial\Omega) := \left\{ g \in H^{1/2}(\partial\Omega) ~\Big|~ \int_{\partial\Omega} g(s) =0\right\}$$
while $ H^1_\diamond(\Omega)$ is the closed subspace of $H^1(\Omega)$ consisting of all functions with zero mean on the boundary, i.e.
$$ H^1_\diamond(\Omega) := \left\{ u \in H^1(\Omega) ~\Big|~ \int_{\partial\Omega} \gamma u =0\right\}.$$

Let us denote by $C^\Omega_\diamond$ the positive constant appearing in the Poincar\'e-Friedrichs inequality (see, for example, \cite{Pechstein})
\begin{align}\label{20-10-15ct1}
C^\Omega_\diamond \int_\Omega \varphi^2 \le \int_\Omega |\nabla \varphi|^2 \mbox{~for all~} \varphi\in  H^1_\diamond(\Omega).
\end{align}
Then for all $q \in \mathcal{Q}$ defined by \eqref{22-3-15ct3},
the coercivity condition
\begin{align}\label{coercivity}
\| \varphi\|^2_{H^1(\Omega)} \le \frac{1+ C^\Omega_\diamond}{C^\Omega_\diamond} \int_\Omega | \nabla \varphi|^2 \le \frac{1+ C^\Omega_\diamond}{C^\Omega_\diamond \underline{q}} \int_\Omega q\nabla\varphi \cdot \nabla\varphi
\end{align}
holds for all $\varphi\in  H^1_\diamond(\Omega)$. Furthermore, since $H^1_0(\Omega) := \left\{ u \in H^1(\Omega) ~\big|~ \gamma u =0\right\} \subset  H^1_\diamond(\Omega)$, the inequality \eqref{coercivity} remains valid for all $\varphi\in H^1_0(\Omega)$.

Finally, for the sake of completeness we briefly introduce the space of functions with bounded total variation;  for more details one may consult \cite{attouch,Giusti}. A scalar function $q\in
L^{1}(\Omega)$ is said to be of bounded total variation if
\begin{align*}
TV(q) := \int_{\Omega}\left|\nabla q\right|
:=\sup\left\{\int_{\Omega} q\divv~ \Xi ~\big|~ \Xi\in C^1_c(\Omega)^d,~
\left| \Xi(x)\right|_{\infty} \le
1,~x\in\Omega\right\}<\infty.
\end{align*}
Here $\left|\cdot\right|_{\infty}$ denotes the $\ell_{\infty}$-norm
on $\mathbb{R}^d$ defined by
$\left|x\right|_{\infty}=\max\limits_{1\le i\le d}\left|x_i\right|$
and $C^1_c(\Omega)$ the space of continuously differentiable functions with compact support in $\Omega$.
The space of all functions in $L^{1}(\Omega)$ with bounded total
variation  is denoted by
$$BV(\Omega)=\left\{ q\in L^{1}(\Omega) ~\Big|~
\int_{\Omega}\left|\nabla q\right|<\infty\right\}$$
which is a Banach space with the norm
$$
\| q \|_{BV(\Omega)}
:= \| q \|_{L^1(\Omega)} + \int_{\Omega} |\nabla
q|.
$$
Furthermore, if $\Omega$ is an open bounded set with
Lipschitz boundary, then $W^{1,1}(\Omega)\varsubsetneq BV(\Omega)$.

\subsection{Neumann operator, Dirichlet operator and Neumann-to-Dirichlet map}\label{operatores}

\subsubsection{Neumann operator} We consider the following Neumann problem
\begin{align}
-\nabla \cdot (q\nabla u) = f \mbox{~in~} \Omega \mbox{~and~} q\nabla u\cdot \vec{n} = j \mbox{~on~} \partial\Omega. \label{22-3-16ct1}
\end{align}
By the coercivity condition \eqref{coercivity} and the Riesz representation theorem, we conclude that for each
$q\in\mathcal{Q}$ and $j\in H^{-1/2}_{-c_f}(\partial\Omega)$
there exists a unique weak solution $u$ of the problem (\ref{22-3-16ct1}) in the sense that $u\in  H^1_\diamond(\Omega)$ and satisfies the identity
\begin{align}
\int_\Omega  q \nabla u \cdot \nabla \varphi =
\left\langle j,\gamma \varphi\right\rangle + (f,\varphi)\label{ct9}
\end{align}
for all $\varphi\in  H^1_\diamond(\Omega)$. By the imposed compatibility condition $\langle j,1 \rangle = -c_f$, i.e.
\begin{equation}\label{compat}
\left\langle j,1\right\rangle + (f,1) =0,
\end{equation}
and the fact that $H^1(\Omega)=H^1_\diamond(\Omega)+\mbox{span}\{1\}$, equation \eqref{ct9} is satisfied for all $\varphi\in  H^1(\Omega)$. Furthermore, this solution satisfies the following estimate
\begin{align}\label{mq5}
\left\|u\right\|_{H^1(\Omega)}
&\le \frac{1+ C^\Omega_\diamond}{C^\Omega_\diamond \underline{q}} \left( \left\|\gamma\right\|_{\mathcal{L}\big(H^1(\Omega),H^{1/2}(\partial\Omega)\big)} \left\|j\right\|_{H^{-1/2}(\partial\Omega)} + \|f\|_{H^{-1}(\Omega)}\right) \notag\\
&\le C_{\mathcal{N}} \left( \left\|j\right\|_{H^{-1/2}(\partial\Omega)} + \|f\|_{H^{-1}(\Omega)}\right),
\end{align}
where
$$ C_{\mathcal{N}} := \frac{1+ C^\Omega_\diamond}{C^\Omega_\diamond \underline{q}} \max \left( 1, \left\|\gamma\right\|_{\mathcal{L}\big(H^1(\Omega),H^{1/2}(\partial\Omega)\big)}\right) .$$
Then for any fixed $j\in H^{-1/2}_{-c_f}(\partial\Omega)$ we can define the {\it Neumann operator}
\begin{align*}
\mathcal{N} : \mathcal{Q}
\rightarrow  H^1_\diamond(\Omega) \mbox{~with~} q\mapsto \mathcal{N}_qj
\end{align*}
which maps each $q \in \mathcal{Q} $ to the unique weak solution $\mathcal{N}_qj := u$ of the problem (\ref{22-3-16ct1}).

\begin{remark} We note that the restriction $j\in H^{-1/2}_{-c_f}(\partial\Omega)$ instead of $j\in H^{-1/2}(\partial\Omega)$ preserves the compatibility condition \eqref{compat} for the pure Neumann problem.
In case this condition fails, the strong form of the problem (\ref{22-3-16ct1}) has no solution. This is the reason why we require $j\in H^{-1/2}_{-c_f}(\partial\Omega)$. However, its weak form, i.e. the variational equation \eqref{ct9}, attains a unique solution independently of the compatibility condition. By working with the weak form only, all results in the present paper remain valid for $j\in H^{-1/2}(\partial\Omega)$.
\end{remark}

\subsubsection{Dirichlet operator}

We now consider the following Dirichlet problem
\begin{align}
-\nabla \cdot (q\nabla v) = f \mbox{~in~} \Omega \mbox{~and~} v = g \mbox{~on~} \partial\Omega.  \label{22-3-16ct2}
\end{align}
For each $q\in\mathcal{Q}$ and $g\in H^{1/2}(\partial\Omega)$, by the coercivity condition \eqref{coercivity}, the problem (\ref{22-3-16ct2}) attains a unique weak solution $v$ in the sense that $v\in H^1(\Omega)$, $\gamma v = g$ and satisfies the identity
\begin{align}
\int_\Omega  q \nabla v \cdot \nabla \psi  =
(f,\psi)  \label{ct9*}
\end{align}
for all $\psi\in H^1_0(\Omega)$. We can rewrite
\begin{align}\label{14-6-17ct3}
v =v_0 +G,
\end{align} 
where $G=\gamma^{-1}g$ and $v_0 \in H^1_0(\Omega)$ is the unique solution to the following variational problem
\begin{align*}
\int_\Omega  q \nabla v_0 \cdot \nabla \psi  = (f,\psi)
 - \int_\Omega  q \nabla G \cdot \nabla \psi
\end{align*}
for all $\psi\in H^1_0(\Omega)$. Since
$$\left\|G\right\|_{H^1(\Omega)} \le \left\|\gamma^{-1}\right\|_{\mathcal{L}\big(H^{1/2}(\partial\Omega),H^1(\Omega)\big)} \left\|g\right\|_{H^{1/2}(\partial\Omega)},$$
we thus obtain the priori estimate
\begin{align}\label{mq5*}
\left\|v\right\|_{H^1(\Omega)} & \le \left\|v_0\right\|_{H^1(\Omega)} + \left\|G\right\|_{H^1(\Omega)}\notag\\
&\le \frac{1+ C^\Omega_\diamond}{C^\Omega_\diamond \underline{q}}\|f\|_{H^{-1}(\Omega)} + \frac{1+ C^\Omega_\diamond}{C^\Omega_\diamond \underline{q}} \overline{q} \left\|G\right\|_{H^1(\Omega)} + \left\|G\right\|_{H^1(\Omega)}\notag\\
&\le C_\mathcal{D} \left( \left\|g\right\|_{H^{1/2}(\partial\Omega)} + \|f\|_{H^{-1}(\Omega)}\right) ,
\end{align}
where
$$C_\mathcal{D} := \max \left( \frac{1+ C^\Omega_\diamond}{C^\Omega_\diamond \underline{q}}, \left(  \frac{1+ C^\Omega_\diamond}{C^\Omega_\diamond \underline{q}} \overline{q} + 1\right) \left\|\gamma^{-1}\right\|_{\mathcal{L}\big(H^{1/2}(\partial\Omega),H^1(\Omega)\big)}\right) . $$
The {\it Dirichlet operator} is for any fixed $g \in H^{1/2}(\partial\Omega)$ defined as
\begin{align*}
\mathcal{D} : \mathcal{Q}
\rightarrow H^1(\Omega) \mbox{~with~} q\mapsto\mathcal{D}_qg
\end{align*}
which maps each $q\in\mathcal{Q}$ to the unique weak solution $\mathcal{D}_qg := v$ of the problem (\ref{22-3-16ct2}).

\subsubsection{Neumann-to-Dirichlet map}

For any fixed $q\in\mathcal{Q}$ we can define the {\it Neumann-to-Dirichlet map}
\begin{align*}
\Lambda_q: H^{-1/2}_{-c_f}(\partial\Omega)
&\rightarrow H^{1/2}_\diamond(\partial\Omega)\\
j&\mapsto \Lambda_qj := \gamma \mathcal{N}_qj.
\end{align*}
Since
\begin{align*}
\int_\Omega  q \nabla \mathcal{N}_q j \cdot \nabla \psi = (f,\psi)
\end{align*}
for all $\psi \in H^1_0(\Omega)$, in view of \eqref{ct9*} we conclude that
$$\Lambda_qj = g \mbox{~if and only if~} \mathcal{N}_qj = \mathcal{D}_qg.$$

\subsection{Identification problem}

The inverse problem is stated as follows.
\begin{center}
{\it Given $f\in H^{-1}(\Omega)$, $\left( j^\dag, g^\dag \right) \in H^{-1/2}_{-c_f}(\partial\Omega) \times H^{1/2}_\diamond(\partial\Omega)$ with $\Lambda_{q^\dag} j^\dag = g^\dag$, find $q^\dag\in \mathcal{Q}$.}
\end{center}
In other words, the problem of interest is, given $f\in H^{-1}(\Omega)$, and a single Neumann-Dirichlet pair $\left( j^\dag, g^\dag \right) \in H^{-1/2}_{-c_f}(\partial\Omega) \times H^{1/2}_\diamond(\partial\Omega)$, to find $q^\dag \in \mathcal{Q} $ and $\Phi^\dag \in  H^1_\diamond(\Omega)$ such that the system \eqref{17-5-16ct1}--\eqref{17-5-16ct3} is satisfied in the weak sense.

\subsection{Total variation regularization}

Assume that $\left(j_\delta,g_\delta\right) \in H^{-1/2}_{-c_f}(\partial\Omega) \times H^{1/2}_\diamond(\partial\Omega)$ is the measured data of the exact boundary values $(j^\dag,g^\dag)$ with
\begin{align}\label{26-3-16ct1}
\big\|j_\delta-j^\dag\big\|_{H^{-1/2}(\partial\Omega)} + \big\|g_\delta-g^\dag\big\|_{H^{1/2}(\partial\Omega)} \le \delta
\end{align}
for some measurement error $\delta >0$.
Our problem is now to reconstruct the conductivity  $q^\dag \in \mathcal{Q}$ from this perturbed data $\left(j_\delta,g_\delta\right)$. For this purpose  we consider the cost functional
\begin{align}\label{26-3-16ct2}
\mathcal{J}_\delta(q) := \int_\Omega q\nabla\left(\mathcal{N}_qj_\delta - \mathcal{D}_qg_\delta\right) \cdot \nabla\left(\mathcal{N}_qj_\delta - \mathcal{D}_qg_\delta\right),
\end{align}
where $\mathcal{N}_qj_\delta$ and $\mathcal{D}_qg_\delta$ is the unique weak solutions of the problems (\ref{22-3-16ct1}) and  (\ref{22-3-16ct2}), respectively, with $j$ in (\ref{22-3-16ct1}) and $g$ in (\ref{22-3-16ct2}) being replaced by $j_\delta$ and $g_\delta$. Furthermore, to estimate the possibly discontinuous conductivity, we here use the total variation regularization (cf., e.g., \cite{BurgerOsher13,Dobson96,EHNBuch}), i.e., we consider the minimization problem
$$\min_{q\in\mathcal{Q}_{ad}} \Upsilon_{\rho,\delta} (q) := \min_{q\in\mathcal{Q}_{ad}} \left( \mathcal{J}_\delta(q) + \rho \int_\Omega | \nabla q|\right), \eqno \left(\mathcal{P}_{\rho,\delta}\right)$$
where
$$\mathcal{Q}_{ad} := \mathcal{Q} \cap BV(\Omega)$$
is the admissible set of the sought conductivities.

\begin{remark}\label{rem:boundarydata}
The noise model \eqref{26-3-16ct1} is to some extent an idealized one, since in practice, measurement precision might be different for the current $j$ and the voltage $g$, and, more importantly, it will be first of all be given with respect to some $L^p$ norm (e.g., $p=2$ corresponding to normally and $p=\infty$ to uniformly distributed noise) rather than in $H^{-1/2}(\partial\Omega) \times H^{1/2}(\partial\Omega)$. While the Neumann data part is unproblematic, by continuity of the embedding of $L^p(\partial\Omega)$ in $H^{-1/2}(\partial\Omega)$ for $p\geq 2\frac{d-1}{d}$, we can obtain an $H^{1/2}(\partial\Omega)$ version of the originally $L^p(\partial\Omega)$ Dirichlet data e.g. by Tikhonov regularization (cf. \cite{EHNBuch} and the references therein) as follows. For simplicity, we restrict ourselves to the Hilbert space case $p=2$ and assume that we have measurements $\tilde{g}_{\delta_g}\in L^2(\partial\Omega)$ such that
$$
\|\tilde{g}_{\delta_g}-g^\dagger\|_{L^2(\partial\Omega)}\leq\delta_g
$$
Tikhonov regularization applied to the embedding operator
$K:H^{1/2}(\partial\Omega)\to L^2(\partial\Omega)$ amounts to finding a minimizer $g_\alpha^{\delta_g}$ of
$$
\min_{g\in H^{1/2}(\partial\Omega)} \|Kg-\tilde{g}_{\delta_g}\|_{L^2(\partial\Omega)}^2
+\alpha \|g\|_{H^{1/2}(\partial\Omega)}^2
$$
where we use
$$
\|g\|_{H^{1/2}(\partial\Omega)}:=\|\gamma^{-1}g\|_{H^1(\Omega)}
=\left(\int_\Omega (|\nabla \gamma^{-1}g|^2+|\gamma^{-1}g|^2)\, dx\right)^{1/2}
$$
as a norm on $H^{1/2}(\partial\Omega)$.
The first order optimality conditions for this quadratic minimization problem yield
$$
\int_{\partial\Omega} \phi(g_\alpha^{\delta_g}-\tilde{g}_{\delta_g})\, ds
+\alpha \int_\Omega (\nabla \gamma^{-1} g_\alpha^{\delta_g}\cdot \nabla\gamma^{-1}\phi + \gamma^{-1} g_\alpha^{\delta_g}\gamma^{-1}\phi)\, dx =0 \mbox{ for all }\phi\in H^{1/2}(\partial\Omega),
$$
which is equivalent to
$$
\int_{\partial\Omega} \gamma\varphi(\gamma w-\tilde{g}_{\delta_g})\, ds
+\alpha \int_\Omega (\nabla w\cdot \nabla\varphi + w\varphi)\, dx =0 \mbox{ for all }\varphi\in H^1(\Omega),
$$
for $w=\gamma^{-1}g_\alpha^{\delta_g}$, i.e., the weak form of the Robin problem
\begin{equation}\label{eqw}
\begin{cases}
-\Delta w+w &= 0 \mbox{~in~} \Omega,\\
 \alpha\nabla w \cdot \vec{n} + w&= \tilde{g}_{\delta_g} \mbox{~on~} \partial\Omega.
\end{cases}
\end{equation}
Thus, according to well-known results from regularization theory (cf., e.g. \cite{EHNBuch}), the smoothed version $g_\delta:=g_\alpha^{\delta_g}=\gamma w$ (where $w$ weakly solves \eqref{eqw}) of $\tilde{g}_{\delta_g}$ converges to $g^\dagger$ as $\delta_g$ tends to zero, provided the regularization parameter $\alpha=\alpha(\delta_g,\tilde{g}_{\delta_g})$ is chosen appropriately. The latter can, e.g., be done by the discrepancy principle, where $\alpha$ is chosen such that
$$
\|Kg_\alpha^{\delta_g}-\tilde{g}_{\delta_g}\|_{L^2(\partial\Omega)}^2
=\int_{\partial\Omega} |g_\alpha^{\delta_g}-\tilde{g}_{\delta_g}|^2\, dx\sim\delta_g^2.
$$
\\
We also wish to mention the complete electrode model cf., e.g., \cite{SomersaloCheneyIsaacson92}, which fully takes into account the fact that current and voltage are typically not measured pointwise on the whole boundary, but via a set of finitely many electrodes with finite geometric extensions as well as contact impedances.
\end{remark}

\subsection{Auxiliary results}

Now we summarize some useful properties of the Neumann and Dirichlet operators. The proof of the following result is based on standard arguments
and therefore omitted.

\begin{lemma} \label{q-diff}
Let $(j,g)\in H^{-1/2}_{-c_f}(\partial\Omega) \times H^{1/2}_\diamond(\partial\Omega)$ be fixed.

(i) The Neumann  operator $\mathcal{N} : \mathcal{Q} \subset L^\infty(\Omega) \rightarrow  H^1_\diamond(\Omega)$ is continuously Fr\'echet differentiable on the set
$\mathcal{Q}$. For each $q\in \mathcal{Q}$ the action of the Fr\'echet derivative in direction $\xi\in L^\infty(\Omega)$ denoted by $\eta_{\mathcal{N}} := \mathcal{N}'_qj(\xi):=\mathcal{N}'(q)\xi$ is the unique weak solution in $ H^1_\diamond(\Omega)$ to the Neumann problem
\begin{align*}
-\nabla \cdot (q\nabla \eta_{\mathcal{N}}) = \nabla \cdot (\xi\nabla \mathcal{N}_qj) \mbox{~in~} \Omega \mbox{~and~} q\nabla \eta_{\mathcal{N}}\cdot \vec{n} = -\xi\nabla \mathcal{N}_qj \cdot\vec{n} \mbox{~on~} \partial\Omega
\end{align*}
in the sense that the identity
\begin{align}
\int_{\Omega}q \nabla \eta_{\mathcal{N}} \cdot \nabla \varphi = -\int_{\Omega} \xi\nabla \mathcal{N}_q j \cdot\nabla \varphi \label{m6}
\end{align}
holds for all $\varphi\in H^1_\diamond(\Omega)$.  Furthermore, the following estimate is fulfilled
\begin{align}
\| \eta_{\mathcal{N}} \|_{H^1(\Omega)}\le \frac{\left( 1+ C^\Omega_\diamond\right) C_{\mathcal{N}}}{C^\Omega_\diamond \underline{q}}
\left( \left\|j\right\|_{H^{-1/2}(\partial\Omega)} + \|f\|_{H^{-1}(\Omega)}\right)
\|\xi\|_{L^\infty(\Omega)}.\label{m6*}
\end{align}

(ii) The Dirichlet  operator $\mathcal{D} : \mathcal{Q} \subset L^\infty(\Omega) \rightarrow  H^1_\diamond(\Omega)$ is continuously Fr\'echet differentiable on the set
$\mathcal{Q}$. For each $q\in \mathcal{Q}$ the action of the Fr\'echet derivative in direction $\xi\in L^\infty(\Omega)$ denoted by $\eta_{\mathcal{D}} := \mathcal{D}'_qg(\xi)=: \mathcal{D}'(q)\xi$ is the unique weak solution in $H^1_0(\Omega)$ to the Dirichlet problem
\begin{align*}
-\nabla \cdot (q\nabla \eta_{\mathcal{D}}) = \nabla \cdot (\xi\nabla \mathcal{D}_qg) \mbox{~in~} \Omega \mbox{~and~} \eta_{\mathcal{D}} =  0 \mbox{~on~} \partial\Omega
\end{align*}
in the sense that it satisfies the equation
\begin{align*}
\int_{\Omega}q \nabla \eta_{\mathcal{D}} \cdot \nabla \psi = -\int_{\Omega} \xi\nabla \mathcal{D}_qg \cdot\nabla \psi
\end{align*}
for all $\psi \in H^1_0(\Omega)$.  Furthermore, the following estimate is
fulfilled
\begin{align*}
\| \eta_{\mathcal{D}} \|_{H^1(\Omega)}\le \frac{\left( 1+ C^\Omega_\diamond\right) C_{\mathcal{D}}}{C^\Omega_\diamond \underline{q}}
\left( \left\|g\right\|_{H^{1/2}(\partial\Omega)} + \|f\|_{H^{-1}(\Omega)}\right)
\|\xi\|_{L^\infty(\Omega)}.
\end{align*}
\end{lemma}

\begin{lemma}\label{weakly conv.}
If the sequence $\left( q_n\right)\subset \mathcal{Q}$
converges to $q$ in the $L^1(\Omega)$-norm, then $q\in\mathcal{Q}$ and for any fixed $(j_\delta,g_\delta)\in H^{-1/2}_{-c_f}(\partial\Omega) \times H^{1/2}_\diamond(\partial\Omega)$ the sequence $\left(\mathcal{N}_{q_n}j_\delta,\mathcal{D}_{q_n}g_\delta\right)$ converges to $\left(\mathcal{N}_{q}j_\delta,\mathcal{D}_{q}g_\delta\right)$ in the $H^1(\Omega)\times H^1(\Omega)$-norm. Furthermore, there holds
\begin{align*}
\limn \mathcal{J}_\delta\left(q_n\right) = \mathcal{J}_\delta\big(q\big),
\end{align*}
where the functional $\mathcal{J}_\delta$ is defined in \eqref{26-3-16ct2}.
\end{lemma}

\begin{proof}
Since $\left(q_n\right)\subset \mathcal{Q}$ converges to $q$ in the
$L^1(\Omega)$-norm, up to a subsequence we assume that it converges to $q$
a.e. in $\Omega$, which implies that $q\in\mathcal{Q}$.
For all $\varphi \in  H^1_\diamond(\Omega)$ we infer from (\ref{ct9}) that
\begin{align*}
\int_\Omega  q_n \nabla\mathcal{N}_{q_n}j_\delta \cdot \nabla \varphi
&= \left\langle j_\delta,\gamma \varphi\right\rangle +(f,\varphi)
%\\&
= \int_\Omega  q \nabla\mathcal{N}_{q}j_\delta \cdot \nabla \varphi
\end{align*}
and so that
\begin{align}\label{6-4-16ct1}
\int_\Omega  q_n \nabla \left(\mathcal{N}_{q_n}j_\delta - \mathcal{N}_{q}j_\delta\right) \cdot \nabla \varphi = \int_\Omega  \left(q -q_n\right) \nabla\mathcal{N}_{q}j_\delta \cdot \nabla \varphi.
\end{align}
Taking $\varphi = \mathcal{N}_{q_n}j_\delta - \mathcal{N}_{q}j_\delta$, by (\ref{coercivity}), we get
\begin{align*}
\frac{C^\Omega_\diamond \underline{q}}{1+ C^\Omega_\diamond}\left\|  \mathcal{N}_{q_n}j_\delta - \mathcal{N}_{q}j_\delta \right \|^2_{H^1(\Omega)}
&\le \int_\Omega  q_n \nabla \left(\mathcal{N}_{q_n}j_\delta - \mathcal{N}_{q}j_\delta\right) \cdot \nabla \left(\mathcal{N}_{q_n}j_\delta - \mathcal{N}_{q}j_\delta\right)\\
&= \int_\Omega  \left(q -q_n\right) \nabla\mathcal{N}_{q}j_\delta \cdot \nabla \left(\mathcal{N}_{q_n}j_\delta - \mathcal{N}_{q}j_\delta\right)\\
&\le \left( \int_\Omega | q - q_n |^2 \left| \nabla \mathcal{N}_qj_\delta \right|^2 \right)^{1/2} \left( \int_\Omega \left|\nabla \left(\mathcal{N}_{q_n}j_\delta - \mathcal{N}_{q}j_\delta\right)\right|^2\right)^{1/2}
\end{align*}
and so that
\begin{align*}
\left\|  \mathcal{N}_{q_n}j_\delta - \mathcal{N}_{q}j_\delta \right \|_{H^1(\Omega)}
\le \frac{1+ C^\Omega_\diamond}{C^\Omega_\diamond \underline{q}} \left( \int_\Omega | q - q_n |^2 \left| \nabla \mathcal{N}_qj_\delta \right|^2 \right)^{1/2}.
\end{align*}
Hence, by the Lebesgue dominated convergence theorem, we deduce from the last inequality that
\begin{align}\label{6-4-16ct2}
\limn \left\|  \mathcal{N}_{q_n}j_\delta - \mathcal{N}_{q}j_\delta \right \|_{H^1(\Omega)}  = 0.
\end{align}
Similarly to \eqref{6-4-16ct1}, we also get
\begin{align*}
\int_\Omega  q_n \nabla \left(\mathcal{D}_{q_n}g_\delta - \mathcal{D}_{q}g_\delta\right) \cdot \nabla \psi = \int_\Omega  \left(q -q_n\right) \nabla\mathcal{D}_{q}g_\delta \cdot \nabla \psi
\end{align*}
for all $\psi\in H^1_0(\Omega)$. Since $\gamma\mathcal{D}_{q_n}g_\delta=\gamma\mathcal{D}_qg_\delta=g_\delta$, taking $\psi=\mathcal{D}_{q_n}g_\delta -\mathcal{D}_{q}g_\delta\in H^1_0(\Omega)$ in the last equation, we also obtain the limit
\begin{align}\label{6-4-16ct3}
\limn \left\|  \mathcal{D}_{q_n}g_\delta - \mathcal{D}_{q}g_\delta \right \|_{H^1(\Omega)}  = 0.
\end{align}
Next, we rewrite the functional $\mathcal{J}_\delta$ as follows
\begin{align}\label{27-3-16ct1}
\mathcal{J}_\delta\left(q_n\right) &= \int_\Omega q_n\nabla \mathcal{N}_{q_n}j_\delta \cdot \nabla \mathcal{N}_{q_n}j_\delta - 2\int_\Omega q_n\nabla \mathcal{N}_{q_n}j_\delta \cdot \nabla \mathcal{D}_{q_n}g_\delta + \int_\Omega q_n\nabla \mathcal{D}_{q_n}g_\delta \cdot \nabla \mathcal{D}_{q_n}g_\delta \notag\\
&= \left\langle j_\delta,\gamma \mathcal{N}_{q_n}j_\delta\right\rangle + \left(f, \mathcal{N}_{q_n}j_\delta\right) - 2\left( \left\langle j_\delta,g_\delta\right\rangle + \left(f, \mathcal{D}_{q_n}g_\delta\right) \right) + \int_\Omega q_n\nabla \mathcal{D}_{q_n}g_\delta \cdot \nabla \mathcal{D}_{q_n}g_\delta
\end{align}
and, by \eqref{6-4-16ct2}--\eqref{6-4-16ct3}, have that
\begin{align}\label{27-3-16ct2}
\left\langle j_\delta,\gamma \mathcal{N}_{q_n}j_\delta\right\rangle + \left(f, \mathcal{N}_{q_n}j_\delta - 2 \mathcal{D}_{q_n}g_\delta\right) \rightarrow \left\langle j_\delta,\gamma \mathcal{N}_{q}j_\delta\right\rangle + \left(f, \mathcal{N}_{q}j_\delta - 2 \mathcal{D}_{q}g_\delta\right)
\end{align}
as $n$ tends to $\infty$. We now consider the difference
\begin{align*}
\int_\Omega q_n\nabla \mathcal{D}_{q_n}g_\delta \cdot \nabla \mathcal{D}_{q_n}g_\delta &- \int_\Omega q\nabla \mathcal{D}_{q}g_\delta \cdot \nabla \mathcal{D}_{q}g_\delta \\
&= \int_\Omega q_n\nabla \left(\mathcal{D}_{q_n}g_\delta - \mathcal{D}_{q}g_\delta\right) \cdot \nabla \left(\mathcal{D}_{q_n}g_\delta + \mathcal{D}_{q}g_\delta\right) - \int_\Omega \left(q - q_n\right)\nabla \mathcal{D}_{q}g_\delta \cdot \nabla \mathcal{D}_{q}g_\delta
\end{align*}
and note that
\begin{align*}
\int_\Omega \left(q - q_n\right)\nabla \mathcal{D}_{q}g_\delta \cdot \nabla \mathcal{D}_{q}g_\delta \rightarrow 0
\end{align*}
as $n$ goes to $\infty$,
by the Lebesgue dominated convergence theorem.
Furthermore, then applying the Cauchy-Schwarz inequality, we also get that
\begin{align*}
\bigg| \int_\Omega q_n&\nabla \left(\mathcal{D}_{q_n}g_\delta - \mathcal{D}_{q}g_\delta\right) \cdot \nabla \left(\mathcal{D}_{q_n}g_\delta + \mathcal{D}_{q}g_\delta\right) \bigg|\\
&\le \overline{q} \left(  \int_\Omega \left| \nabla \left(\mathcal{D}_{q_n}g_\delta - \mathcal{D}_{q}g_\delta\right) \right|^2\right)^{1/2} \left(  \int_\Omega \left| \nabla \left(\mathcal{D}_{q_n}g_\delta + \mathcal{D}_{q}g_\delta\right) \right|^2\right)^{1/2}\\
&\le \overline{q} \left\| \mathcal{D}_{q_n}g_\delta - \mathcal{D}_{q}g_\delta \right\|_{H^1(\Omega)} \left( \left\| \mathcal{D}_{q_n}g_\delta \right\|_{H^1(\Omega)} + \left\|\mathcal{D}_{q}g_\delta \right\|_{H^1(\Omega)}\right)\rightarrow 0
\end{align*}
as $n$ approaches $\infty$, here we used \eqref{mq5*} and \eqref{6-4-16ct3}. We thus obtain that
\begin{align}\label{27-3-16ct3}
\int_\Omega q_n\nabla \mathcal{D}_{q_n}g_\delta \cdot \nabla \mathcal{D}_{q_n}g_\delta \rightarrow \int_\Omega q\nabla \mathcal{D}_{q}g_\delta \cdot \nabla \mathcal{D}_{q}g_\delta
\end{align}
as $n$ tends to $\infty$. Then we deduce from (\ref{27-3-16ct1})--(\ref{27-3-16ct3}) that
\begin{align*}
\limn \mathcal{J}_\delta\left(q_n\right) &= \left\langle j_\delta,\gamma \mathcal{N}_{q}j_\delta\right\rangle + \left(f, \mathcal{N}_{q}j_\delta\right) -2 \left\langle j_\delta,g_\delta\right\rangle -2 \left(f, \mathcal{D}_{q}g_\delta\right) + \int_\Omega q\nabla \mathcal{D}_{q}g_\delta \cdot \nabla \mathcal{D}_{q}g_\delta\\
&= \int_\Omega q\nabla \mathcal{N}_{q}j_\delta \cdot \nabla \mathcal{N}_{q}j_\delta - 2\int_\Omega q\nabla \mathcal{N}_{q}j_\delta \cdot \nabla \mathcal{D}_{q}g_\delta + \int_\Omega q\nabla \mathcal{D}_{q}g_\delta \cdot \nabla \mathcal{D}_{q}g_\delta\\
&= \mathcal{J}_\delta\big(q\big),
\end{align*}
which finishes the proof.
\end{proof}

\begin{lemma}[\cite{Giusti}]\label{bv1}

(i) Let $\left(q_n\right)$ be a bounded sequence in the $BV(\Omega)$-norm.
Then a subsequence which is denoted by the same symbol and an element
$q\in BV(\Omega)$ exist such that $\left(q_n\right)$ converges to $q$ in the
$L^1(\Omega)$-norm.

(ii) Let $\left(q_n\right)$ be a sequence in $BV(\Omega)$ converging to
$q$ in  the $L^1(\Omega)$-norm. Then $q \in BV(\Omega)$ and
\begin{align}\label{27-7-16ct1}
\int_\Omega \left|\nabla q\right| \le \liminfn \int_\Omega|\nabla q_n|.
\end{align}
\end{lemma}
We mention that equality need not be achieved in \eqref{27-7-16ct1}. Here is a counterexample from \cite{Giusti}. Let $\Omega=(0,2\pi)$ and $q_n(x)=\frac{1}{n}\sin nx$ for $x\in\Omega$ and $n\in\mathbb{N}$. Then $\|q_n\|_{L^1(\Omega)} \to 0$ as $n\to \infty$, but $\int_\Omega|\nabla q_n| =4$ for each $n\in\mathbb{N}$.

Let us quote the following useful result on approximation of $BV$-functions by smooth functions.

\begin{lemma}[\cite{bartels,ChenZou}]\label{appr.BV}
Assume that $w\in BV(\Omega)$. Then for all $\alpha >0$ an element $w^\alpha \in C^\infty(\Omega)$ exists such that
$$\int_\Omega|w-w^\alpha| \le \alpha\int_\Omega|\nabla w|,~ \int_\Omega|\nabla w^\alpha| \le (1+C\alpha)\int_\Omega|\nabla w| \mbox{~and~} \int_\Omega|D^2 w^\alpha| \le C\alpha^{-1}\int_\Omega|\nabla w|,$$
where the positive constant $C$ is independent of $\alpha$.
\end{lemma}

Now, we are in a position to prove the main result of this section

\begin{theorem}\label{existance}
The problem $\left(\mathcal{P}_{\rho,\delta}\right)$ attains a solution $q_{\rho,\delta}$, which is called the regularized solution of the identification problem.
\end{theorem}

\begin{proof}
Let $\left(q_n\right) \subset \mathcal{Q}_{ad}$ be a minimizing sequence of the problem $\left(\mathcal{P}_{\rho,\delta}\right)$, i.e.,
\begin{align}\label{26-3-16ct4}
\limn \left(\mathcal{J}_\delta\left(q_n\right) + \rho \int_\Omega \left| \nabla q_n\right|\right) = \inf_{q\in \mathcal{Q}_{ad}} \left(\mathcal{J}_\delta(q) + \rho \int_\Omega \left| \nabla q\right|\right).
\end{align}
Then, due to Lemma \ref{bv1}, a subsequence which is not relabelled and an element $q \in \mathcal{Q}_{ad}$ exist such that $\left(q_n\right)$ converges to $q$ in the
$L^1(\Omega)$-norm and
\begin{align}\label{26-3-16ct5}
\int_\Omega |\nabla q| \le \liminfn \int_\Omega|\nabla q_n|.
\end{align}
Using Lemma \ref{weakly conv.} and by \eqref{26-3-16ct4}--\eqref{26-3-16ct5}, we obtain that
\begin{align*}
\mathcal{J}_\delta (q) + \rho\int_\Omega |\nabla q|
&\le \limn \mathcal{J}_\delta\left(q_n\right) + \liminfn \rho\int_\Omega|\nabla q_n| \\
&= \liminfn \left(\mathcal{J}_\delta\left(q_n\right) + \rho \int_\Omega \left| \nabla q_n\right|\right) \\
&= \inf_{q\in \mathcal{Q}_{ad}} \left(\mathcal{J}_\delta(q) + \rho \int_\Omega \left| \nabla q\right|\right).
\end{align*}
This means that $q$ is a solution of the problem $\left(\mathcal{P}_{\rho,\delta}\right)$, which finishes the proof.
\end{proof}

\section{Finite element method for the identification problem}\label{Finite element method}

Let $\left(\mathcal{T}^h\right)_{0<h<1}$ be a family of regular and
quasi-uniform triangulations of the domain $\overline{\Omega}$ with the mesh size $h$ such that each vertex of the polygonal boundary $\partial\Omega$ is a node of $\mathcal{T}^h$. For the definition of the discretization space of the state
functions let us denote
\begin{equation*}
\mathcal{V}_1^h := \left\{v^h\in C(\overline\Omega)
~\big|~{v^h}_{|T} \in \mathcal{P}_1(T), ~~\forall
T\in \mathcal{T}^h\right\}
\end{equation*}
and
$$\mathcal{V}_{1,\diamond}^h := \mathcal{V}_1^h \cap  H^1_\diamond(\Omega) \mbox{~and~} \mathcal{V}_{1,0}^h := \mathcal{V}_1^h \cap H^1_0(\Omega) \subset \mathcal{V}_{1,\diamond}^h,$$
where $\mathcal{P}_1$ consists of all polynomial functions of degree
less than or equal to 1.

To go further, we introduce the following modified Cl\'ement's interpolation operator, see \cite{Clement}.

\begin{lemma}\label{moli.data}
An interpolation operator $\Pi^h_\diamond: L^1(\Omega) \rightarrow \mathcal{V}^h_{1,\diamond}$ exists such that
$$\Pi^h_\diamond\varphi^h = \varphi^h \mbox{~for all~} \varphi^h \in \mathcal{V}^h_{1,\diamond}  \mbox{~and~} \Pi^h_\diamond \big(H^1_0(\Omega)\big) \subset \mathcal{V}^h_{1,0} \subset \mathcal{V}^h_{1,\diamond}.$$
Furthermore, it satisfies the properties
\begin{equation}\label{23/10:ct2*}
\lim_{h\to 0} \big\| \vartheta - \Pi^h_\diamond \vartheta
\big\|_{H^1(\Omega)} =0 \enskip \mbox{~for all~} \vartheta \in H^1_\diamond(\Omega)
\end{equation}
and
\begin{equation}\label{23/5:ct1*}
\big\| \vartheta - \Pi^h_\diamond \vartheta \big\|_{H^1(\Omega)} \le
Ch \| \vartheta\|_{H^2(\Omega)} \mbox{~for all~} \vartheta \in H^1_\diamond(\Omega)\cap H^2(\Omega)
\end{equation}
with the positive constant $C$ being independent of $h$ and $\vartheta$.
\end{lemma}

\begin{proof}
It is well known (see \cite{Clement} and some generalizations \cite{Bernardi1,Bernardi2,scott_zhang}) that there is an interpolation operator
$$\Pi^h: L^1(\Omega) \rightarrow \mathcal{V}^h_1 \mbox{~with~} \Pi^h\varphi^h = \varphi^h \mbox{~for all~} \varphi^h \in \mathcal{V}^h_{1}  \mbox{~and~} \Pi^h \big(H^1_0(\Omega)\big) \subset \mathcal{V}^h_{1,0}$$
which satisfies the following properties
\begin{equation}\label{23/10:ct2}
\lim_{h\to 0} \big\| \vartheta - \Pi^h \vartheta
\big\|_{H^1(\Omega)} =0 \enskip \mbox{~for all~} \vartheta \in H^1(\Omega)
\end{equation}
and
\begin{equation}\label{23/5:ct1}
\big\| \vartheta - \Pi^h \vartheta \big\|_{H^1(\Omega)} \le
Ch \| \vartheta\|_{H^2(\Omega)} \mbox{~for all~} \vartheta \in H^2(\Omega).
\end{equation}
We then define for each $\vartheta\in L^1(\Omega)$ 
$$\Pi^h_\diamond\vartheta := \Pi^h\vartheta -\frac{1}{|\partial\Omega|}\int_{\partial\Omega} \gamma \Pi^h\vartheta \in \mathcal{V}^h_{1,\diamond}.$$
Then $\Pi^h_\diamond \big(L^1(\Omega)\big) \subset \mathcal{V}^h_{1,\diamond}$, $\Pi^h_\diamond \varphi^h = \varphi^h$ for all $\varphi^h \in \mathcal{V}^h_{1,\diamond}$  and $\Pi^h_\diamond \big(H^1_0(\Omega)\big) \subset \mathcal{V}^h_{1,0}$. Furthermore, since $\nabla\Pi^h_\diamond\vartheta = \nabla\Pi^h\vartheta$ for all $\vartheta\in L^1(\Omega)$, the properties \eqref{23/10:ct2*}, \eqref{23/5:ct1*} are deduced from
\eqref{23/10:ct2}, \eqref{23/5:ct1}, respectively. The proof is completed.
\end{proof}
We remark that the operator $\Pi^h$ in the above proof satisfies the estimate $\|\vartheta-\Pi^h\vartheta\|_{H^k(\Omega)} \le Ch^{l-k}\|\vartheta\|_{H^l(\Omega)}$ for $0\le k\le l\le 2$ and $\vartheta\in H^l(\Omega)$ (see \cite{Clement}) which implies that 
\begin{align}\label{21-6-17ct1}
\left\|\Pi^h_\diamond\vartheta\right\|_{H^1(\Omega)} \le C\left\|\vartheta\right\|_{H^1(\Omega)} \quad\mbox{for all}\quad \vartheta\in H^1(\Omega),
\end{align}
an estimate that is required for the proof of part (ii) of the following proposition.

Similarly to the continuous case we have the following result.
\begin{proposition}\label{14-6-17ct2}
(i) Let $q$ be in $\mathcal{Q}$ and $j$ be in $H^{-1/2}_{-c_f}(\partial\Omega)$. Then the variational equation
\begin{align}
\int_\Omega q\nabla u^h \cdot \nabla \varphi^h = \left\langle j,\gamma\varphi^h\right\rangle + \left( f,\varphi^h\right) \mbox{~for all~} \varphi^h\in
\mathcal{V}_{1,\diamond}^h \label{10/4:ct1}
\end{align}
admits a unique solution $u^h\in \mathcal{V}_{1,\diamond}^h$. Furthermore, there holds
\begin{align}
\big\|u^h\big\|_{H^1(\Omega)}\le C_{\mathcal{N}}
\left( \left\|j\right\|_{H^{-1/2}(\partial\Omega)} + \|f\|_{H^{-1}(\Omega)}\right). \label{18/5:ct1}
\end{align}
(ii) Let $q$ be in $\mathcal{Q}$ and $g$ be in $H^{1/2}_\diamond(\partial\Omega)$. Then the equation
\begin{align}
\int_\Omega q\nabla v^h \cdot \nabla \psi^h = \left( f,\psi^h\right) \mbox{~for all~} \psi^h\in
\mathcal{V}_{1,0}^h \label{10/4:ct1*}
\end{align}
with $\gamma v^h = \gamma\big(\Pi^h_\diamond (\gamma^{-1}g)\big)$ has a unique solution $v^h\in \mathcal{V}_{1,\diamond}^h$. Furthermore, the stability estimate
\begin{align}
\big\|v^h\big\|_{H^1(\Omega)}\le \bar{C}_{\mathcal{D}}
\left( \left\|g\right\|_{H^{1/2}(\partial\Omega)} + \|f\|_{H^{-1}(\Omega)}\right) \label{18/5:ct1*}
\end{align}
is satisfied, where $\bar{C}_\mathcal{D} := \max \left( \frac{1+ C^\Omega_\diamond}{C^\Omega_\diamond \underline{q}}, \left(  \frac{1+ C^\Omega_\diamond}{C^\Omega_\diamond \underline{q}} \overline{q} + 1\right) \left\|\Pi^h_\diamond\right\|_{\mathcal{L}\big(H^1(\Omega),H^1(\Omega)\big)} \left\|\gamma^{-1}\right\|_{\mathcal{L}\big(H^{1/2}(\partial\Omega),H^1(\Omega)\big)}\right)$.
\end{proposition}

Let $u$ and $u^h$ be solutions to \eqref{ct9} and \eqref{10/4:ct1}, respectively. Due to the standard theory of the finite element method (see, for example, \cite{Brenner_Scott,Ciarlet}), the estimate
\begin{align}\label{12-6-17ct1*}
\|u-u^h\|_{H^1(\Omega)} \le Ch\|u\|_{H^2(\Omega)} 
\end{align}
holds in case $u\in H^2(\Omega)$, where the positive constant $C$ is independent of $h$ and $u$.

Assume that $v$ and $v^h$ are the solutions to \eqref{ct9*} and \eqref{10/4:ct1*}, where $v\in H^2(\Omega)$, we then have (see, for example, \cite[Section 5.4]{Brenner_Scott}) that
\begin{align*}
\|v-v^h\|_{H^1(\Omega)} \le \inf_{\psi^h\in
\mathcal{V}_{1,0}^h}\|v-\gamma^{-1}g-\psi^h\|_{H^1(\Omega)} + 2\|\gamma^{-1}g - \Pi^h_\diamond (\gamma^{-1}g)\|_{H^1(\Omega)}.
\end{align*}
Since $v\in H^2(\Omega)$, it follows that $g=\gamma v\in H^{3/2}(\Omega)$ and so $\gamma^{-1}g \in H^2(\Omega)$. Due to the approximation property of the finite dimensional spaces $\mathcal{V}_{1,0}^h \subset H^1_0(\Omega)$ (which states that $\inf_{\psi^h\in
\mathcal{V}_{1,0}^h}\|\psi-\psi^h\|_{H^1(\Omega)} \le C h\|\psi\|_{H^2(\Omega)}$ for each $\psi\in H^2(\Omega)\cap H^1_0(\Omega)$, where the constant $C$ is independent of $h$ and $\psi$) and \eqref{23/5:ct1}, we deduce 
\begin{align}\label{12-6-17ct1}
\|v-v^h\|_{H^1(\Omega)} \le Ch\left(\|v\|_{H^2(\Omega)} + \|\gamma^{-1} g\|_{H^2(\Omega)}\right).
\end{align}
We also mention that above we approximate the Dirichlet boundary condition $g$ by $g^h := \gamma\big(\Pi^h_\diamond (\gamma^{-1}g)$. There exist some different choices for the approximation $g^h$; for example, the $L^2$-projection of $g$ on the set $\mathcal{S}^h_{\partial\Omega} := \{\gamma \varphi^h ~|~ \varphi^h \in \mathcal{V}^h_1\}$, or the Lagrange interpolation of $g$ in $\mathcal{S}^h_{\partial\Omega}$ in case $g$ being smooth enough (see \cite{fix83} for more details).

\begin{definition}\label{discrete_solution}
(i) For any fixed $j\in H^{-1/2}_{-c_f}(\partial\Omega)$ the operator $\mathcal{N}^h: \mathcal{Q} \rightarrow
\mathcal{V}_{1,\diamond}^h$ mapping each $q \in  \mathcal{Q}$
to the unique solution $u^h =: \mathcal{N}^h_qj$ of the variational equation \eqref{10/4:ct1}
is called the  {\it discrete Neumann operator}.

(ii) For any fixed $g\in H^{1/2}_\diamond(\partial\Omega)$ the operator $\mathcal{D}^h: \mathcal{Q} \rightarrow
\mathcal{V}_{1,\diamond}^h$ mapping each $q \in  \mathcal{Q}$
to the unique solution $v^h =: \mathcal{D}^h_q g$ of the variational equation \eqref{10/4:ct1*}
is called the  {\it discrete Dirichlet operator}.
\end{definition}

Next, the discretization space for the sought conductivity is defined by
$$\mathcal{Q}^h_{ad} := \mathcal{Q}\cap \mathcal{V}^h_1 \subset \mathcal{Q} \cap BV(\Omega) = \mathcal{Q}_{ad}.$$
Then, using the discrete operators $\mathcal{N}^h$ and $\mathcal{D}^h$ in Definition \ref{discrete_solution}, we introduce
the discrete cost functional
%corresponding to $\Upsilon_{\rho,\delta}$
\begin{align}\label{27-3-16ct4}
\Upsilon^h_{\rho,\delta} (q):= \mathcal{J}_\delta^h (q) + \rho \int_\Omega \sqrt{\left|\nabla q\right|^2+\epsilon^h},
\end{align}
where $q \in \mathcal{Q}^h_{ad}$, $\epsilon^h$ is a positive function of the mesh size $h$ satisfying $\lim_{h\to0} \epsilon^h =0$ and
\begin{equation}\label{29/6:ct9}
\mathcal{J}_\delta^h (q):= \int_\Omega q\nabla \left(\mathcal{N}^h_qj_\delta-
\mathcal{D}^h_qg_\delta\right) \cdot \nabla \left(\mathcal{N}^h_qj_\delta-
\mathcal{D}^h_qg_\delta\right) \mbox{~with~} q\in\mathcal{Q}.
\end{equation}
The positive function $\epsilon^h$ above acts as a smoothing parameter for the total variation.

\begin{theorem}\label{dis.solution}
For any fixed $h$, $\rho$ and $\delta$ the minimization problem
$$\min_{q\in\mathcal{Q}^h_{ad}} \Upsilon^h_{\rho,\delta} (q) \eqno\left(\mathcal{P}^h_{\rho,\delta}\right)$$
attains a solution $q^h_{\rho,\delta}$, which is called the discrete regularized solution of the identiﬁcation problem.
\end{theorem}

\begin{proof}
We first note that $\mathcal{Q}^h_{ad}$ is a compact subset of the finite dimensional space $\mathcal{V}^h_1$. Let $\left(q_n\right) \subset \mathcal{Q}^h_{ad}$ be a minimizing sequence of the problem $\left(\mathcal{P}^h_{\rho,\delta}\right)$, i.e.,
\begin{align}\label{26-3-16ct4*}
\limn \Upsilon^h_{\rho,\delta} \left(q_n\right) = \inf_{q\in \mathcal{Q}^h_{ad}} \Upsilon^h_{\rho,\delta} (q).
\end{align}
Then a subsequence of $\left(q_n\right)$ which is denoted by the same symbol and an element $q \in \mathcal{Q}^h_{ad}$ exist such that $\left(q_n\right)$ converges to $q$ in the $H^1(\Omega)$-norm. We have that
\begin{align}\label{26-3-16ct5*}
&\bigg| \int_\Omega \sqrt{\left|\nabla q_n\right|^2+\epsilon^h} - \int_\Omega \sqrt{\left|\nabla q\right|^2+\epsilon^h} \bigg|
\le \int_\Omega\frac{\left| \left|\nabla q_n\right|^2 - \left|\nabla q\right|^2\right|}{ \sqrt{\left|\nabla q_n\right|^2 +\epsilon^h} + \sqrt{\left|\nabla q\right|^2 +\epsilon^h}} \notag\\
&\le \frac{1}{2\sqrt{\epsilon^h}}\int_\Omega\left| \nabla q_n - \nabla q\right| \left| \nabla q_n + \nabla q\right|
\le \frac{1}{2\sqrt{\epsilon^h}} \left(\int_\Omega\left| \nabla q_n - \nabla q\right|^2\right)^{1/2} \left(\int_\Omega\left| \nabla q_n + \nabla q\right|^2\right)^{1/2} \notag\\
&\le \frac{1}{2\sqrt{\epsilon^h}} \left\| q_n - q\right\|_{H^1(\Omega)} \left(\left\| q_n-q\right\|_{H^1(\Omega)} + 2\left\| q\right\|_{H^1(\Omega)}\right)
\rightarrow 0 \mbox{~as~} n\to\infty.
\end{align}
On the other hand, similarly to Lemma \ref{weakly conv.}, we can prove that the sequence $\left(\mathcal{N}^h_{q_n}j_\delta,\mathcal{D}^h_{q_n}g_\delta\right)$ converges to $\left(\mathcal{N}^h_{q}j_\delta,\mathcal{D}^h_{q}g_\delta\right)$ in the $H^1(\Omega)\times H^1(\Omega)$-norm as $n$ goes to $\infty$ and then obtain
\begin{align}\label{27-3-16ct5}
\limn \mathcal{J}^h_\delta\left(q_n\right) = \mathcal{J}^h_\delta(q).
\end{align}
Thus, it follows from \eqref{26-3-16ct4*}--\eqref{27-3-16ct5} that
\begin{align*}
\Upsilon^h_{\rho,\delta} (q)
= \limn \Upsilon^h_{\rho,\delta} \left(q_n\right)
= \inf_{q\in \mathcal{Q}^h_{ad}} \Upsilon^h_{\rho,\delta} (q),
\end{align*}
which finishes the proof.
\end{proof}

\section{Convergence}\label{Stability}

From now on $C$ is a generic positive constant which is independent of the mesh size $h$ of $\mathcal{T}^h$, the noise level $\delta$ and the regularization parameter $\rho$. The following result shows the stability of the finite element method for the regularized identification problem.

\begin{theorem}\label{odinh1}
Let $(h_n)_n$ be a sequence with $\limn h_n = 0$ and $\left( j_{\delta_n}, g_{\delta_n}\right)$ be a sequence in $H^{-1/2}_{-c_f}(\partial\Omega) \times H^{1/2}_\diamond(\partial\Omega)$ converging to $\left( j_\delta, g_\delta\right) $ in the  $H^{-1/2}(\partial\Omega) \times H^{1/2}(\partial\Omega)$-norm. For a fixed regularization parameter $\rho >0$ let
$q^{h_n}_{\rho,\delta_n} \in \mathcal{Q}^{h_n}_{ad}$ be a minimizer of
$\left(\mathcal{P}^{h_n}_{\rho,\delta_n} \right)$ for each $n\in\mathbb{N}$.
Then a subsequence of $\big(q^{h_n}_{\rho,\delta_n} \big)$ not relabelled and an element $q_{\rho,\delta} \in \mathcal{Q}_{ad}$ exist such that
$$\limn \big\| q^{h_n}_{\rho,\delta_n} - q_{\rho,\delta}\big\|_{L^1(\Omega)} = 0 \mbox{~and~} \limn \int_\Omega \big|\nabla q^{h_n}_{\rho,\delta_n} \big| = \int_\Omega \left| \nabla q_{\rho,\delta} \right|.$$
Furthermore, $q_{\rho,\delta}$ is a solution to
$\left( \mathcal{P}_{\rho,\delta} \right)$.
\end{theorem}

To prove the theorem, we need the auxiliary results, starting with the following estimates.

\begin{lemma}\label{6-5-16ct5}
Let $(j_1,g_1)$ and $(j_2,g_2)$ be arbitrary in $  H^{-1/2}_{-c_f}(\partial\Omega) \times H^{1/2}_\diamond(\partial\Omega)$. Then the estimates
\begin{align}\label{7-5-16ct1}
\left\|\mathcal{N}^{h}_{q} j_1 - \mathcal{N}^{h}_{q} j_2 \right\|_{H^1(\Omega)} \le \frac{1+ C^\Omega_\diamond}{C^\Omega_\diamond \underline{q}} \left\|\gamma\right\|_{\mathcal{L}\big(H^1(\Omega),H^{1/2}(\partial\Omega)\big)} \left\|j_1 -j_2 \right\|_{H^{-1/2}(\partial\Omega)}
\end{align}
and
\begin{align}\label{7-5-16ct2}
\left\|\mathcal{D}^{h}_{q} g_1 - \mathcal{D}^{h}_{q} g_2 \right\|_{H^1(\Omega)} \le\left(  \frac{1+ C^\Omega_\diamond}{C^\Omega_\diamond \underline{q}} \overline{q} + 1\right) \left\|\Pi^h_\diamond\right\|_{\mathcal{L}\big(H^1(\Omega),H^1(\Omega)\big)}\left\|\gamma^{-1}\right\|_{\mathcal{L}\big(H^{1/2}(\partial\Omega),H^1(\Omega)\big)} \left\|g_1 -g_2 \right\|_{H^{1/2}(\partial\Omega)}
\end{align}
hold for all $q\in \mathcal{Q}$ and $h>0$.
\end{lemma}

\begin{proof}
According to the definition of the discrete Neumann operator, we have for all $\varphi^{h}\in
\mathcal{V}_{1,\diamond}^{h}$ that
$$\int_\Omega  q \nabla \mathcal{N}^{h}_{q} j_i \cdot \nabla \varphi^{h} =
\left\langle j_i,\gamma \varphi^{h}\right\rangle + \left( f,\varphi^h\right) \mbox{~with~} i=1,2.$$
Thus, $\Phi^{h}_{\mathcal{N}} := \mathcal{N}^{h}_{q} j_1 - \mathcal{N}^{h}_{q} j_2$ is the unique solution to the variational problem
$$\int_\Omega  q \nabla \Phi^{h}_{\mathcal{N}} \cdot \nabla \varphi^{h} =
\left\langle j_1-j_2,\gamma \varphi^{h}\right\rangle$$
for all $\varphi^{h}\in
\mathcal{V}_{1,\diamond}^{h}$
and so that \eqref{7-5-16ct1} follows. Similarly, we also obtain \eqref{7-5-16ct2}, which finishes the proof.
\end{proof}

\begin{lemma}\label{weakly conv.dis}
Let $(h_n)_n$ be a sequence with $\limn h_n = 0$ and $\left( j_{\delta_n}, g_{\delta_n}\right) \subset H^{-1/2}_{-c_f}(\partial\Omega) \times H^{1/2}_\diamond(\partial\Omega)$ be a sequence converging to $\left( j_\delta, g_\delta\right) $ in the  $H^{-1/2}(\partial\Omega) \times H^{1/2}(\partial\Omega)$-norm. Then for any fixed $q\in\mathcal{Q}$ the limit
\begin{align}\label{20-5-16ct2}
\limn \mathcal{J}^{h_n}_{\delta_n} (q) = \mathcal{J}_{\delta} (q)
\end{align}
holds. Furthermore, if $(q_n)$ is a sequence in $\mathcal{Q}$ which converges to $q$ in the $L^1(\Omega)$-norm, then the sequence $\left(\mathcal{N}^{h_n}_{q_n}j_{\delta_n},\mathcal{D}^{h_n}_{q_n}g_{\delta_n}\right)$ converges to $\left(\mathcal{N}_{q}j_\delta,\mathcal{D}_{q}g_\delta\right)$ in the $H^1(\Omega)\times H^1(\Omega)$-norm and the limit
\begin{align} \label{20-5-16ct3}
\limn \mathcal{J}^{h_n}_{\delta_n} (q_n) = \mathcal{J}_{\delta} (q)
\end{align}
also holds.
\end{lemma}

\begin{proof}
We get for any fixed $q\in\mathcal{Q}$ that
$$\mathcal{N}^{h_n}_{q} j_{\delta_n} - \mathcal{D}^{h_n}_{q} g_{\delta_n} = \left( \mathcal{N}_{q} j_{\delta} - \mathcal{D}_{q} g_{\delta}\right) + \left( \mathcal{N}^{h_n}_{q} j_{\delta_n} - \mathcal{N}_{q} j_{\delta} + \mathcal{D}_{q} g_{\delta} - \mathcal{D}^{h_n}_{q} g_{\delta_n} \right).$$
Thus, with $\Phi_n := \mathcal{N}^{h_n}_{q} j_{\delta_n} - \mathcal{N}_{q} j_{\delta} + \mathcal{D}_{q} g_{\delta} - \mathcal{D}^{h_n}_{q} g_{\delta_n}$ we have
\begin{align*}
\mathcal{J}^{h_n}_{\delta_n} (q) &= \int_\Omega q\nabla\left(\mathcal{N}^{h_n}_{q} j_{\delta_n} - \mathcal{D}^{h_n}_{q} g_{\delta_n}\right) \cdot \nabla\left(\mathcal{N}^{h_n}_{q} j_{\delta_n} - \mathcal{D}^{h_n}_{q} g_{\delta_n}\right) \\
&= \mathcal{J}_{\delta} (q) +\int_\Omega q\nabla \Phi_n \cdot \nabla \Phi_n  +2 \int_\Omega q\nabla \left( \mathcal{N}_{q} j_{\delta} - \mathcal{D}_{q} g_{\delta}\right) \cdot \nabla \Phi_n.
\end{align*}
Applying Lemma \ref{6-5-16ct5}, we infer that
\begin{align*}
\left\|\mathcal{N}^{h_n}_{q} j_{\delta_n} - \mathcal{N}_{q} j_{\delta} \right\|_{H^1(\Omega)} &\le \left\|\mathcal{N}^{h_n}_{q} j_{\delta} - \mathcal{N}_{q} j_{\delta} \right\|_{H^1(\Omega)} + \left\|\mathcal{N}^{h_n}_{q} j_{\delta_n} - \mathcal{N}^{h_n}_{q} j_{\delta} \right\|_{H^1(\Omega)} \\
&\le \left\|\mathcal{N}^{h_n}_{q} j_{\delta} - \mathcal{N}_{q} j_{\delta} \right\|_{H^1(\Omega)} + C \left\|j_{\delta_n} - j_{\delta} \right\|_{H^{-1/2}(\partial\Omega)} \rightarrow 0 \mbox{~as~} n\to\infty,
\end{align*}
where we used the limit
$$\limn \left\|\mathcal{N}^{h_n}_{q} j_{\delta} - \mathcal{N}_{q} j_{\delta} \right\|_{H^1(\Omega)} =0,$$
due to the standard theory (see, for example, \cite{Brenner_Scott,Ciarlet}). Similarly, we also have
\begin{align*}
\left\|\mathcal{D}^{h_n}_{q} g_{\delta_n} - \mathcal{D}_{q} g_{\delta} \right\|_{H^1(\Omega)} \rightarrow 0 \mbox{~as~} n\to\infty.
\end{align*}
We thus get that
\begin{align*}
\|\Phi_n\|_{H^1(\Omega)} &\le \left\|\mathcal{N}^{h_n}_{q} j_{\delta_n} - \mathcal{N}_{q} j_{\delta} \right\|_{H^1(\Omega)} + \left\|\mathcal{D}_{q} g_{\delta} - \mathcal{D}^{h_n}_{q} g_{\delta_n} \right\|_{H^1(\Omega)}  \rightarrow 0 \mbox{~as~} n\to\infty.
\end{align*}
Therefore, we obtain that
\begin{align*}
\limn \left| \int_\Omega q\nabla \Phi_n \cdot \nabla \Phi_n  +2 \int_\Omega q\nabla \left( \mathcal{N}_{q} j_{\delta} - \mathcal{D}_{q} g_{\delta}\right) \cdot \nabla \Phi_n \right| &\le C \limn \left( \|\Phi_n\|^2_{H^1(\Omega)} + \|\Phi_n\|_{H^1(\Omega)}\right)=0
\end{align*}
and \eqref{20-5-16ct2} then follows.

Next, for $q_n$ converging to $q$ in $L^1(\Omega)$, hence, along a subsequence again denoted by $(q_n)_n$, pointwise almost everywhere, by (\ref{10/4:ct1}) and \eqref{ct9}, we have
\begin{align*}
\int_\Omega  q_n \nabla\mathcal{N}^{h_n}_{q_n}j_{\delta_n} \cdot \nabla \varphi^{h_n}
&= \big\langle j_{\delta_n},\gamma \varphi^{h_n}\big\rangle + \big(f, \varphi^{h_n}\big)= \big\langle j_{\delta},\gamma \varphi^{h_n}\big\rangle + \big(f, \varphi^{h_n}\big) + \big\langle j_{\delta_n} - j_\delta,\gamma \varphi^{h_n}\big\rangle\\
&= \int_\Omega  q \nabla\mathcal{N}_{q}j_\delta \cdot \nabla \varphi^{h_n} + \big\langle j_{\delta_n} - j_\delta,\gamma \varphi^{h_n}\big\rangle
\end{align*}
for all $\varphi^{h_n}\in \mathcal{V}_{1,\diamond}^{h_n}$ which implies that
\begin{align}\label{decomplemlim}
\int_\Omega  q_n \nabla &\left(\mathcal{N}^{h_n}_{q_n}j_{\delta_n} -  \Pi_\diamond^{h_n}\mathcal{N}_{q}j_\delta\right) \cdot \nabla \varphi^{h_n} \notag \\
&= \int_\Omega  q \nabla\mathcal{N}_{q}j_\delta \cdot \nabla \varphi^{h_n} - \int_\Omega  q_n \nabla \Pi_\diamond^{h_n} \mathcal{N}_{q}j_\delta \cdot \nabla \varphi^{h_n} + \big\langle j_{\delta_n} - j_\delta,\gamma \varphi^{h_n}\big\rangle \notag\\
&= \int_\Omega  \big(q-q_n\big) \nabla\mathcal{N}_{q}j_\delta \cdot \nabla \varphi^{h_n} + \int_\Omega  q_n \nabla\left( \mathcal{N}_{q}j_\delta - \Pi_\diamond^{h_n} \mathcal{N}_{q}j_\delta \right)\cdot \nabla \varphi^{h_n} + \big\langle j_{\delta_n} - j_\delta,\gamma \varphi^{h_n}\big\rangle,
\end{align}
where the operator $\Pi_\diamond^{h_n}$ is defined according to Lemma \ref{moli.data}. Taking $\varphi^{h_n} = \mathcal{N}^{h_n}_{q_n}j_{\delta_n} - \Pi_\diamond^{h_n}\mathcal{N}_{q}j_\delta \in \mathcal{V}_{1,\diamond}^h$, by (\ref{coercivity}) and using the Cauchy-Schwarz inequality, we get
\begin{align*}
\frac{C^\Omega_\diamond \underline{q}}{1+ C^\Omega_\diamond} & \left\| \mathcal{N}^{h_n}_{q_n}j_{\delta_n} - \Pi_\diamond^{h_n}\mathcal{N}_{q}j_\delta \right\|^2_{H^1(\Omega)}\\
&\le \left( \int_\Omega | q - q_n |^2 \left| \nabla \mathcal{N}_qj_\delta \right|^2 \right)^{1/2} \left\| \mathcal{N}^{h_n}_{q_n}j_{\delta_n} - \Pi_\diamond^{h_n}\mathcal{N}_{q}j_\delta \right\|_{H^1(\Omega)} \\
&~\quad + \overline{q} \left\| \mathcal{N}_{q}j_\delta - \Pi_\diamond^{h_n} \mathcal{N}_{q}j_\delta \right\|_{H^1(\Omega)} \left\| \mathcal{N}^{h_n}_{q_n}j_{\delta_n} - \Pi_\diamond^{h_n}\mathcal{N}_{q}j_\delta \right\|_{H^1(\Omega)}\\
&~\quad + \left\|\gamma\right\|_{\mathcal{L}\big(H^1(\Omega),H^{1/2}(\partial\Omega)\big)} \left\|j_{\delta_n} - j_\delta\right\|_{H^{-1/2}(\partial\Omega)} \left\| \mathcal{N}^{h_n}_{q_n}j_{\delta_n} - \Pi_\diamond^{h_n}\mathcal{N}_{q}j_\delta \right\|_{H^1(\Omega)}
\end{align*}
and so that
\begin{align*}
\frac{C^\Omega_\diamond \underline{q}}{1+ C^\Omega_\diamond}\Big\| \mathcal{N}^{h_n}_{q_n}j_{\delta_n} - \Pi_\diamond^{h_n}\mathcal{N}_{q}j_\delta \Big\|_{H^1(\Omega)}
&\le\left( \int_\Omega | q - q_n |^2 \left| \nabla \mathcal{N}_qj_\delta \right|^2 \right)^{1/2}  + \overline{q} \left\| \mathcal{N}_{q}j_\delta - \Pi_\diamond^{h_n} \mathcal{N}_{q}j_\delta \right\|_{H^1(\Omega)} \notag\\
&~\quad+ \left\|\gamma\right\|_{\mathcal{L}\big(H^1(\Omega),H^{1/2}(\partial\Omega)\big)} \left\|j_{\delta_n} - j_\delta\right\|_{H^{-1/2}(\partial\Omega)}\rightarrow 0 \mbox{~as~}n\to\infty,
\end{align*}
by the Lebesgue dominated convergence theorem and \eqref{23/10:ct2*}. Thus, we infer from the triangle inequality that
\begin{align*}
\left\| \mathcal{N}^{h_n}_{q_n}j_{\delta_n} - \mathcal{N}_{q}j_\delta \right\|_{H^1(\Omega)} &\le \left\| \mathcal{N}^{h_n}_{q_n}j_{\delta_n} - \Pi_\diamond^{h_n}\mathcal{N}_{q}j_\delta \right\|_{H^1(\Omega)} + \left\| \Pi_\diamond^{h_n} \mathcal{N}_{q}j_\delta -\mathcal{N}_{q}j_\delta \right\|_{H^1(\Omega)}\rightarrow 0 \mbox{~as~} n\to \infty.
\end{align*}
Similarly, using \eqref{ct9*} and \eqref{10/4:ct1*}, for all $\psi^{h_n}\in \mathcal{V}_{1,0}^{h_n}$ we arrive at
\begin{align}\label{19-6-17ct1}
\int_\Omega  q_n \nabla &\left(\mathcal{D}^{h_n}_{q_n}g_{\delta} -  \Pi_\diamond^{h_n}\mathcal{D}_{q}g_{\delta}\right) \cdot \nabla \psi^{h_n}= \int_\Omega  \big(q-q_n\big) \nabla\mathcal{D}_{q}g_{\delta} \cdot \nabla \psi^{h_n} + \int_\Omega  q_n \nabla\left( \mathcal{D}_{q}g_{\delta} - \Pi_\diamond^{h_n} \mathcal{D}_{q}g_{\delta} \right)\cdot \nabla \psi^{h_n}.
\end{align}
We have 
\begin{align}\label{14-6-17ct1}
\gamma\mathcal{D}^{h_n}_{q_n}g_{\delta}=\gamma \big(\Pi^{h_n}_\diamond (\gamma^{-1}g_\delta)\big),
\end{align}
by Proposition \ref{14-6-17ct2} (ii). On the other hand, in view of \eqref{14-6-17ct3}, we get $\mathcal{D}_qg_\delta =v_0 +\gamma^{-1}g_\delta$ with $v_0\in H^1_0(\Omega)$, and therefore
\begin{align}\label{14-6-17ct5}
\gamma\big(\Pi_\diamond^{h_n} \mathcal{D}_{q}g_{\delta}\big) = \gamma\big(\Pi_\diamond^{h_n} (v_0 + \gamma^{-1} g_{\delta})\big) = \gamma\big(\Pi_\diamond^{h_n} v_0\big) + \gamma\big(\Pi_\diamond^{h_n} (\gamma^{-1} g_{\delta})\big) = \gamma\big(\Pi_\diamond^{h_n} (\gamma^{-1} g_{\delta})\big),
\end{align}
since $\gamma\big(\Pi_\diamond^{h_n} v_0\big)=0$.
It follows from \eqref{14-6-17ct1}--\eqref{14-6-17ct5} that 
$$\psi^{h_n}_* : = \mathcal{D}^{h_n}_{q_n}g_{\delta} - \Pi_\diamond^{h_n} \mathcal{D}_{q}g_{\delta} \in \mathcal{V}_{1,0}^{h_n}$$
Taking $\psi^{h_n} := \psi^{h_n}_*$ in the above equation \eqref{19-6-17ct1}, it is deduced that
$$\limn \left\| \mathcal{D}^{h_n}_{q_n}g_{\delta} - \Pi_\diamond^{h_n}\mathcal{D}_{q}g_{\delta} \right\|_{H^1(\Omega)} = 0.$$
Using  Lemma \ref{6-5-16ct5}, we therefore obtain that
\begin{align*}
\Big\| \mathcal{D}^{h_n}_{q_n}g_{\delta_n} &- \mathcal{D}_{q}g_{\delta} \Big\|_{H^1(\Omega)} \\
&\le \left\| \mathcal{D}^{h_n}_{q_n}g_{\delta_n} - \mathcal{D}^{h_n}_{q_n}g_{\delta} \right\|_{H^1(\Omega)}  + \left\| \mathcal{D}^{h_n}_{q_n}g_{\delta} - \Pi_\diamond^{h_n}\mathcal{D}_{q}g_{\delta} \right\|_{H^1(\Omega)} + \left\| \Pi_\diamond^{h_n}\mathcal{D}_{q}g_{\delta} - \mathcal{D}_{q}g_{\delta} \right\|_{H^1(\Omega)}\\
&\le C\left\|g_{\delta_n} - g_\delta\right\|_{H^{1/2}(\partial\Omega)} + \left\| \mathcal{D}^{h_n}_{q_n}g_{\delta} - \Pi_\diamond^{h_n}\mathcal{D}_{q}g_{\delta} \right\|_{H^1(\Omega)} + \left\| \Pi_\diamond^{h_n}\mathcal{D}_{q}g_{\delta} - \mathcal{D}_{q}g_{\delta} \right\|_{H^1(\Omega)}\\
&\rightarrow 0 \mbox{~as~} n\to \infty.
\end{align*}
Since $\big( q_n\big)$
converges to $q$ in the $L^1(\Omega)$-norm while the sequence $\left(\mathcal{N}^{h_n}_{q_n}j_{\delta_n},\mathcal{D}^{h_n}_{q_n}g_{\delta_n}\right)$ converges to $\left(\mathcal{N}_{q}j_\delta,\mathcal{D}_{q}g_\delta\right)$ in the $H^1(\Omega)\times H^1(\Omega)$-norm, we conclude, similarly to the proof of Lemma \ref{weakly conv.} that
\begin{align*}
\limn \mathcal{J}_{\delta_n}^{h_n} \big(q_n\big)
= \mathcal{J}_\delta(q),
\end{align*}
which finishes the proof.
\end{proof}

\begin{proof}[Proof of Theorem \ref{odinh1}]
To simplify notation we write $q_n := q^{h_n}_{\rho,\delta_n}$. Let $q \in\mathcal{Q}_{ad}$ be arbitrary. Using Lemma \ref{appr.BV}, for any fixed $\alpha\in(0,1)$ an element  $q^\alpha \in C^\infty(\Omega)$ exists such that
\begin{align}\label{28-3-18ct1}
\left\|q-q^\alpha\right\|_{L^1(\Omega)} \le \overline{C}\alpha \mbox{~and~} \int_\Omega \left|\nabla q^\alpha\right| \le \overline{C}\alpha+\int_\Omega \left|\nabla q\right|,
\end{align}
where the positive constant $\overline{C}$ is independent of $\alpha$. Setting
$$q^\alpha_P := \max\left( \underline{q}, \min\left(q^\alpha, \overline{q} \right) \right) \in W^{1,\infty}(\Omega)\cap\mathcal{Q} \subset \mathcal{Q}_{ad}\mbox{~and~} q^\alpha_n := I^{h_n}_1 q^\alpha_P \in \mathcal{Q}^{h_n}_{ad},$$
where
$$I^h_1 : W^{1,p}(\Omega)\hookrightarrow C(\overline{\Omega}) \to \mathcal{V}^h_1 \mbox{~with~} p>d$$
is the usual nodal value interpolation operator. Since the sequence $\big(q^\alpha_n\big)$ converges to $q^\alpha_P$ in the $H^1(\Omega)$-norm as $n$ tends to $\infty$ (see, for example, \cite{Brenner_Scott,Ciarlet}), we get the equation
\begin{align}\label{22-6-16ct2}
\limn \int_\Omega \sqrt{\left|\nabla q^\alpha_n\right|^2+\epsilon^{h_n}} = \int_\Omega |\nabla q^\alpha_P|.
\end{align}
Indeed, we have that
\begin{align*}
\left| \int_\Omega \sqrt{\left|\nabla q^\alpha_n\right|^2+\epsilon^{h_n}} - \int_\Omega \left|\nabla q^\alpha_n\right| \right|
\le \int_\Omega \frac{\epsilon^{h_n}}{\sqrt{\left|\nabla q^\alpha_n\right|^2+\epsilon^{h_n}} + \left|\nabla q^\alpha_n\right|}
\le |\Omega|\sqrt{\epsilon^{h_n}} \ \rightarrow 0 \mbox{~as~} n\to \infty
\end{align*}
and by the reverse triangle as well as the Cauchy Schwarz inequality
\begin{align*}
\left| \int_\Omega \left|\nabla q^\alpha_n\right| - \int_\Omega |\nabla q^\alpha_P| \right|
&\le  \left\|\nabla q^\alpha_n - \nabla q^\alpha_P\right\|_{L^1(\Omega)}
\le  |\Omega|^{1/2} \left\|\nabla q^\alpha_n - \nabla q^\alpha_P\right\|_{L^2(\Omega)}
\rightarrow 0 \mbox{~as~} n\to \infty
\end{align*}
so that \eqref{22-6-16ct2} follows from the triangle inequality.
By \eqref{28-3-18ct1} and the fact that $q^\alpha_P$ is constant on $\{x\in\Omega ~|~ q^\alpha_P(x) \not= q^\alpha(x)\}$, we have that
\begin{align}\label{6-4-16ct5}
\int_\Omega |\nabla q^\alpha_P|
&= \int_{ \{x\in\Omega ~|~ q^\alpha_P(x) = q^\alpha(x)\}} |\nabla q^\alpha_P| \le \int_\Omega |\nabla q^\alpha| \le \overline{C}\alpha + \int_\Omega |\nabla q|.
\end{align}
By the optimality of $q_n$, we get for all $n\in \mathbb{N}$ that
\begin{align}\label{23-6-16ct1}
\mathcal{J}^{h_n}_{\delta_n} \left(q_n\right) + \rho \int_{\Omega} \sqrt{|\nabla q_n|^2 + \epsilon^{h_n}}
\le \mathcal{J}^{h_n}_{\delta_n} \left(q^\alpha_n\right) + \rho \int_{\Omega} \sqrt{|\nabla q^\alpha_n|^2 + \epsilon^{h_n}},
\end{align}
where, by \eqref{18/5:ct1} and \eqref{18/5:ct1*},
\begin{align*}
\mathcal{J}^{h_n}_{\delta_n} \left(q^\alpha_n\right) \le C
\end{align*}
holds for some $C$ independent of $n$ and $\alpha$. We then deduce from \eqref{22-6-16ct2}--\eqref{23-6-16ct1} that
$$\int_{\Omega} |\nabla q_n| \le \int_{\Omega} \sqrt{|\nabla q_n|^2 + \epsilon^{h_n}}\le C(\rho)$$
for another constant $C(\rho)$ independent of $n$ and $\alpha$, but depending on $\rho$, so the sequence $\left(q_n\right)$ is bounded in the $BV(\Omega)$-norm. Thus, by Lemma \ref{bv1}, a subsequence which is denoted by the same symbol and an element $\widehat{q}\in \mathcal{Q}_{ad}$ exist such that $\left(q_n\right)$ converges to $\widehat{q}$ in the
$L^1(\Omega)$-norm and
\begin{align}\label{29-3-16ct1}
\int_\Omega \left|\nabla \widehat{q}\right| \le \liminfn \int_\Omega\left|\nabla q_n\right| \le \liminfn\int_\Omega \sqrt{\left|\nabla q_n\right|^2+\epsilon^{h_n}}.
\end{align}
Furthermore, due to Lemma \ref{weakly conv.dis} we get that
\begin{align}\label{29-3-16ct2*}
\mathcal{J}_\delta(\widehat{q}) = \limn \mathcal{J}^{h_n}_{\delta_n}\left(q_n\right)
\end{align}
and
\begin{align}\label{29-3-16ct2}
\mathcal{J}_\delta \left( q^\alpha_P \right) = \limn \mathcal{J}^{h_n}_{\delta_n}\left(q^\alpha_n\right).
\end{align}
Therefore, by \eqref{22-6-16ct2}--\eqref{29-3-16ct2}, we have that
\begin{align}\label{22-6-16ct3}
\mathcal{J}_\delta(\widehat{q}) +\rho\int_\Omega \left|\nabla \widehat{q}\right| &\le \limn \mathcal{J}^{h_n}_{\delta_n}\left(q_n\right) + \liminfn \rho\int_\Omega \sqrt{\left|\nabla q_n\right|^2+\epsilon^{h_n}}, \mbox{~by \eqref{29-3-16ct2*} and \eqref{29-3-16ct1}~} \notag\\
&=\liminfn \left(\mathcal{J}^{h_n}_{\delta_n}\left(q_n\right) + \rho\int_\Omega \sqrt{\left|\nabla q_n\right|^2+\epsilon^{h_n}}\right)\notag\\
&\le \liminfn \left(  \mathcal{J}^{h_n}_{\delta_n} \left(q^\alpha_n\right)  + \rho \int_{\Omega} \sqrt{|\nabla q^\alpha_n|^2 + \epsilon^{h_n}}\right), \mbox{~by \eqref{23-6-16ct1}} \notag\\
&= \mathcal{J}_\delta \left( q^\alpha_P \right) + \rho\int_\Omega |\nabla q^\alpha_P|, \mbox{~by \eqref{29-3-16ct2} and \eqref{22-6-16ct2}~} \notag\\
&\le \mathcal{J}_\delta \left( q^\alpha_P \right) + \rho \int_\Omega |\nabla q | +\overline{C}\alpha\rho, \mbox{~by \eqref{6-4-16ct5}}.
\end{align}
Now, by the definition of $q^\alpha_P$, we get $|q^\alpha_P -q| \le |q^\alpha - q|$ a.e. in $\Omega$ and therefore
$$\left\|q^\alpha_P - q\right\|_{L^1(\Omega)} \le \left\|q^\alpha - q\right\|_{L^1(\Omega)} \le \overline{C}\alpha.$$
Sending $\alpha$ to zero in the last inequality and applying Lemma \ref{weakly conv.},
we arrive at
\begin{align*}
\mathcal{J}_\delta(\widehat{q}) +\rho\int_\Omega \left|\nabla \widehat{q}\right|
\le \mathcal{J}_\delta(q) +\rho\int_\Omega \left|\nabla q\right| ,
\end{align*}
where $\widehat{q}\in \mathcal{Q}_{ad}$ and $q\in\mathcal{Q}_{ad}$ is arbitrary. This means that $\widehat{q}$ is a solution to
$\left( \mathcal{P}_{\rho,\delta} \right)$ and $\left(q_n\right)$ converges to $\widehat{q}$ in the
$L^1(\Omega)$-norm.

Next, as above, from $\widehat{q}$ we can obtain $\widehat{q}^\alpha$, $\widehat{q}^\alpha_P$, $\widehat{q}^\alpha_n$ and note that $\big(\widehat{q}^\alpha_n\big)$ converges to $\widehat{q}^\alpha_P$ in the $H^1(\Omega)$-norm, so also in the $L^1(\Omega)$-norm, as $n$ tends to $\infty$ while $\big(\widehat{q}^\alpha_P\big)$ converges to $\widehat{q}$ in the $L^1(\Omega)$-norm as $\alpha$ tends to $0$. Then,  by the optimality of $q_n$, we have that
\begin{align}\label{23-6-16ct2}
\mathcal{J}^{h_n}_{\delta_n} \left(q_n\right) + \rho \int_{\Omega} \sqrt{|\nabla q_n|^2 + \epsilon^{h_n}}
&\le \mathcal{J}^{h_n}_{\delta_n} \left(\widehat{q}^\alpha_n\right) + \rho \int_{\Omega} \sqrt{|\nabla \widehat{q}^\alpha_n|^2 + \epsilon^{h_n}}.
\end{align}
By \eqref{29-3-16ct2*}, we then obtain that
\begin{align*}
\rho\limsupn \int_\Omega\left|\nabla q_n\right|
&= \limn \mathcal{J}^{h_n}_{\delta_n}\left(q_n\right) + \rho\limsupn \int_\Omega\left|\nabla q_n\right| - \mathcal{J}_\delta(\widehat{q})\\
&\le \limsupn \left( \mathcal{J}^{h_n}_{\delta_n}\left(q_n\right) + \rho \int_\Omega \sqrt{\left|\nabla q_n\right|^2 + \epsilon^{h_n}} \right)  - \mathcal{J}_\delta(\widehat{q}) \\
&\le \limsupn \left( \mathcal{J}^{h_n}_{\delta_n} \left(\widehat{q}^\alpha_n\right) + \rho \int_{\Omega} \sqrt{|\nabla \widehat{q}^\alpha_n|^2 + \epsilon^{h_n}}\right)   - \mathcal{J}_\delta(\widehat{q}), \mbox{~by \eqref{23-6-16ct2}}\\
&= \mathcal{J}_\delta(\widehat{q}^\alpha_P) + \rho \int_\Omega |\nabla \widehat{q}^\alpha_P| - \mathcal{J}_\delta(\widehat{q}), \mbox{~ by Lemma \ref{weakly conv.dis}}\\
&\le \mathcal{J}_\delta(\widehat{q}^\alpha_P) + \rho \int_\Omega |\nabla \widehat{q}| +\overline{C}\alpha\rho - \mathcal{J}_\delta(\widehat{q}).
\end{align*}
Sending $\alpha$ to zero, we obtain from the last inequality that $\limsupn \int_\Omega\left|\nabla q_n\right| \le \int_\Omega |\nabla \widehat{q}|$. Combining this with \eqref{29-3-16ct1}, we conclude $\limn \int_\Omega\left|\nabla q_n\right| = \int_\Omega |\nabla \widehat{q}|$, which finishes the proof.
\end{proof}

Next we show convergence of the regularized finite element approximations to a solution of the identification problem. Before doing so, we introduce the notion of the {\it total variation-minimizing solution}.

\begin{lemma}\label{identification_problem}
The problem
$$
\min_{q \in \mathcal{I}_{\mathcal{Q}_{ad}} \left(j^\dag,g^\dag\right)} \int_{\Omega}\left|\nabla
q\right| \eqno\left(\mathcal{IP}\right)
$$
attains a solution, which is called the {\it total
variation-minimizing solution} of the identification problem,
where
\begin{align}\label{3-5-16ct1}
\mathcal{I}_{\mathcal{Q}_{ad}} \left(j^\dag,g^\dag\right) :=\left\{q\in \mathcal{Q}_{ad} ~\big|~ \Lambda_qj^\dag = g^\dag \right\} = \left\{q\in \mathcal{Q}_{ad} ~\big|~ \mathcal{N}_qj^\dag = \mathcal{D}_qg^\dag \right\}.
\end{align}
\end{lemma}

\begin{proof}
By our assumption on consistency of the exact boundary data, the set $\mathcal{I}_{\mathcal{Q}_{ad}} \left(j^\dag,g^\dag\right)$ is non-empty. Let $\left(q_n\right) \subset \mathcal{I}_{\mathcal{Q}_{ad}} \left(j^\dag,g^\dag\right)$ be a minimizing sequence of the problem $\left(\mathcal{IP}\right)$, i.e.,
\begin{align}\label{26-3-16ct3*}
\limn \int_{\Omega}\left|\nabla
q_n\right| = \inf_{q\in \mathcal{I}_{\mathcal{Q}_{ad}} \left(j^\dag,g^\dag\right)} \int_{\Omega}\left|\nabla
q\right|.
\end{align}
Then due to Lemma \ref{bv1}, a subsequence which is denoted by the same symbol and an element
$\tilde{q}\in \mathcal{Q}_{ad}$ exist such that $\left(q_n\right)$ converges to $\tilde{q}$ in the
$L^1(\Omega)$-norm and
\begin{align}\label{26-3-16ct3}
\int_\Omega \left|\nabla \tilde{q}\right| \le \limn \int_\Omega|\nabla q_n|.
\end{align}
On the other hand, by Lemma \ref{weakly conv.}, we have that
$$\left(\mathcal{N}_{q_n}j^\dag,\mathcal{D}_{q_n}g^\dag\right) \rightarrow \left(\mathcal{N}_{\tilde{q}}j^\dag,\mathcal{D}_{\tilde{q}}g^\dag\right) \mbox{~in the~} H^1(\Omega)\times H^1(\Omega)\mbox{-norm}.$$
By the definition of the set $\mathcal{I}_{\mathcal{Q}_{ad}} \left(j^\dag,g^\dag\right)$, we get that $\mathcal{N}_{q_n}j^\dag=\mathcal{D}_{q_n}g^\dag$ which implies $\mathcal{N}_{\tilde{q}}j^\dag=\mathcal{D}_{\tilde{q}}g^\dag$.
Combining this with \eqref{26-3-16ct3*} and \eqref{26-3-16ct3}, we conclude that
\begin{align*}
\int_{\Omega}\left|\nabla \tilde{q}\right| \le \inf_{q\in \mathcal{I}_{\mathcal{Q}_{ad}} \left(j^\dag,g^\dag\right)} \int_{\Omega}\left|\nabla q\right|,
\end{align*}
where $\tilde{q}\in\mathcal{I}_{\mathcal{Q}_{ad}} \left(j^\dag,g^\dag\right)$, which finishes the proof.
\end{proof}

\begin{remark}
Note that due to the lack of strict convexity of the cost functional and the admissible  set, a solution of $\left(\mathcal{IP}\right)$ may be nonunique.
\end{remark}

\begin{lemma}\label{23-6-16ct5}
For any fixed $q\in \mathcal{Q}_{ad}$ an element $\widehat{q}^h \in \mathcal{Q}^h_{ad}$ exists such that
\begin{align}\label{22-7-16ct5}
\big\|\widehat{q}^h-q\big\|_{L^1(\Omega)} \le Ch|\log h|
\end{align}
and
\begin{align}\label{22-7-16ct6}
\lim_{h\to 0}\int_\Omega \big|\nabla \widehat{q}^h\big| = \int_\Omega |\nabla q|.
\end{align}
In case $q\in W^{1,p}(\Omega)\hookrightarrow C(\overline{\Omega})$ with $p>d$ the above element $\widehat{q}^h$ can be taken as $I^h_1q$.
\end{lemma}

\begin{proof}
According to Lemma \ref{appr.BV}, for any fixed $\alpha\in(0,1)$ an element  $q^\alpha \in C^\infty(\Omega)$ exists such that
\begin{align*}
\left\|q-q^\alpha\right\|_{L^1(\Omega)} \le \overline{C}\alpha,~ \int_\Omega \left|\nabla q^\alpha\right| \le \overline{C}\alpha+\int_\Omega \left|\nabla q\right| \mbox{~and~} \int_\Omega|D^2 q^\alpha| \le \overline{C}\alpha^{-1}\int_\Omega|\nabla q|,
\end{align*}
where the positive constant $\overline{C}$ is independent of $\alpha$.
Setting
$$q^\alpha_P := \max\left( \underline{q}, \min\left(q^\alpha, \overline{q} \right) \right) \in W^{1,\infty}(\Omega)\cap\mathcal{Q} \subset \mathcal{Q}_{ad}\mbox{~and~} \widehat{q}^h := I^h_1 q^\alpha_P \in \mathcal{Q}^h_{ad},$$
we then have
$$\big| \widehat{q}^h(x) - q(x)\big| = \big| I^h_1q^\alpha(x) - q(x)\big| \mbox{~a.e. in~} \Omega_1 := \{x\in\Omega ~|~ \underline{q} \le q^\alpha \le \overline{q}\}$$
and
$$\big| \widehat{q}^h(x) - q(x)\big| \le \left| q^\alpha(x) - q(x)\right| \mbox{~a.e. in~} \Omega\setminus \Omega_1.$$
We thus have, using for example \cite[Theorem 4.4.20]{Brenner_Scott}, with an another positive constant $C$ independent of $\alpha$ that
\begin{align*}
\big\|\widehat{q}^h-q\big\|_{L^1(\Omega)}
&\le \big\|I^h_1 q^\alpha -q\big\|_{L^1(\Omega_1)} +\left\|q^\alpha-q\right\|_{L^1(\Omega\setminus\Omega_1)}\notag\\
&\le \big\|I^h_1 q^\alpha -q^\alpha\big\|_{L^1(\Omega)} + \left\|q-q^\alpha\right\|_{L^1(\Omega_1)}+\left\|q^\alpha-q\right\|_{L^1(\Omega\setminus\Omega_1)} \notag\\
&\le Ch \int_\Omega \left|\nabla q^\alpha\right| + \left\|q-q^\alpha\right\|_{L^1(\Omega)}\notag\\
&\le Ch \left( \overline{C}\alpha + \int_\Omega \left|\nabla q\right|\right) +\overline{C}\alpha
\le C(h +\alpha)\\
&\le Ch|\log h|
\end{align*}
for $\alpha\sim h|\log h|$. To establish the limit \eqref{22-7-16ct6} we first note that
\begin{align}\label{22-7-16ct1}
\int_\Omega \big| \nabla I^h_1 q^\alpha_P\big| \le \int_\Omega \big| \nabla I^h_1 q^\alpha\big|.
\end{align}
Indeed, we rewrite
\begin{align}\label{22-7-16ct2}
\int_\Omega \big| \nabla I^h_1 q^\alpha_P\big| = \sum_{T\in\mathcal{T}^h_1}\int_T \big| \nabla I^h_1 q^\alpha_P\big| + \sum_{T\in\mathcal{T}^h_2}\int_T \big| \nabla I^h_1 q^\alpha_P\big| + \sum_{T\in\mathcal{T}^h \setminus \big(\mathcal{T}^h_1 \cup \mathcal{T}^h_2\big)}\int_T \big| \nabla I^h_1 q^\alpha_P\big|,
\end{align}
where
$\mathcal{T}^h_1$ includes all triangles $T\in\mathcal{T}^h$ with its vertices $x_1, \ldots, x_d,x_{d+1}$ at which
$$ \mbox{either~} q^\alpha(x_1), \ldots, q^\alpha(x_{d+1}) < \underline{q} \mbox{~or~} q^\alpha(x_1), \ldots, q^\alpha(x_{d+1}) > \overline{q}$$
while $\mathcal{T}^h_2$ consists all triangles $T\in\mathcal{T}^h$ with its vertices $x_1, \ldots, x_d,x_{d+1}$ at which
$$ q^\alpha(x_1), \ldots, q^\alpha(x_{d+1}) \in [ \underline{q}, \overline{q}].$$
We then have that
\begin{align}\label{22-7-16ct3}
\sum_{T\in\mathcal{T}^h_1}\int_T \big| \nabla I^h_1 q^\alpha_P\big| =0 \mbox{~and~} \sum_{T\in\mathcal{T}^h_2}\int_T \big| \nabla I^h_1 q^\alpha_P\big| = \sum_{T\in\mathcal{T}^h_2}\int_T \big| \nabla I^h_1 q^\alpha\big|.
\end{align}
Now let $T\in\mathcal{T}^h \setminus \big(\mathcal{T}^h_1 \cup \mathcal{T}^h_2\big)$ be arbitrary. In Cartesian coordinate system $Oxz$ with $x\in \R^d$ we consider plane surfaces $z=I^h_1 q^\alpha_P(x)$ and $z=I^h_1 q^\alpha(x)$ with $x\in T$ and denote by $\vec{m}_P$ and $\vec{m}$ the constant unit normal on these surfaces in the upward $z$ direction, respectively. By the definition of the projection $q^\alpha_P$, we get $0 < \widehat{(Oz,\vec{m}_P)} \le \widehat{(Oz,\vec{m})} < \pi/2$
and so that
$0< \cos \widehat{(Oz,\vec{m})} \le \cos \widehat{(Oz,\vec{m}_P)} <1.$
Since
$$\cos \widehat{(Oz,\vec{m})} = \frac{1}{\sqrt{\big| \nabla I^h_1 q^\alpha\big|^2 + 1}} \mbox{~and~} \cos \widehat{(Oz,\vec{m}_P)} = \frac{1}{\sqrt{\big| \nabla I^h_1 q^\alpha_P\big|^2 + 1}},$$
it follows that $\big| \nabla I^h_1 q^\alpha(x)\big| \ge \big| \nabla I^h_1 q^\alpha_P(x)\big|$ for all $x\in T$. We thus have that
\begin{align}\label{22-7-16ct4}
\sum_{T\in\mathcal{T}^h \setminus \big(\mathcal{T}^h_1 \cup \mathcal{T}^h_2\big)}\int_T \big| \nabla I^h_1 q^\alpha_P\big| \le \sum_{T\in\mathcal{T}^h \setminus \big(\mathcal{T}^h_1 \cup \mathcal{T}^h_2\big)}\int_T \big| \nabla I^h_1 q^\alpha\big|.
\end{align}
The inequality \eqref{22-7-16ct1} is then directly deduced from \eqref{22-7-16ct2}--\eqref{22-7-16ct4}. We therefore have with a constant $C$ independent of $\alpha$ that
\begin{align*}
\int_\Omega \big| \nabla \widehat{q}^h\big| - \int_\Omega \big| \nabla q\big| &= \int_\Omega \big| \nabla I^h_1 q^\alpha_P\big| - \int_\Omega \big| \nabla q\big| \le \int_\Omega \big| \nabla I^h_1 q^\alpha\big| - \int_\Omega \big| \nabla q\big| \\
&\le \int_\Omega \big| \nabla \big( I^h_1 q^\alpha -q^\alpha\big) \big| + \int_\Omega \big| \nabla q^\alpha \big|- \int_\Omega \big| \nabla q\big| \\
&\le Ch \int_\Omega \big| D^2 q^\alpha \big| + \overline{C}\alpha \\
&\le C\overline{C}h\alpha^{-1} \int_\Omega |\nabla q| + \overline{C}\alpha\\
&\le C\left( |\log h|^{-1} + h|\log h|\right)  \to 0 \mbox{~as~} h\to 0 \mbox{~and for~} \alpha\sim h|\log h|.
\end{align*}
Combining this with \eqref{22-7-16ct5} and Lemma \ref{bv1}, we obtain that
$$\int_\Omega | \nabla q| \le \liminf_{h\to 0} \int_\Omega \big| \nabla \widehat{q}^h\big| \le \limsup_{h\to 0} \int_\Omega \big| \nabla \widehat{q}^h\big|\le \int_\Omega | \nabla q|,$$
which finishes the proof.
\end{proof}

\begin{lemma}\label{7-4-16ct1}
Let $(q,j,g) \in \mathcal{Q}_{ad} \times H^{-1/2}_{-c_f}(\partial\Omega) \times H^{1/2}_\diamond(\partial\Omega)$ be arbitrary. Then the convergence
\begin{align*}
\widehat{\varrho}^h_{q} \left(j,g\right) := \left\| \mathcal{N}^h_{\widehat{q}^h}j - \mathcal{N}_q j \right\|_{H^1(\Omega)} + \left\| \mathcal{D}^h_{\widehat{q}^h}g - \mathcal{D}_q g \right\|_{H^1(\Omega)} \rightarrow 0 \mbox{~as~} h\to 0
\end{align*}
holds, where $ \widehat{q}^h$ is generated from $q$ according to Lemma \ref{23-6-16ct5}.
\end{lemma}

\begin{proof}
The assertion follows directly from Lemma \ref{weakly conv.dis} and Lemma \ref{23-6-16ct5}.
\end{proof}

Additional smoothness assumptions enable an error estimate of $\widehat{\varrho}^h_{q} \left(j,g\right)$.

\begin{lemma}\label{7-4-16ct3}
Let $(q,j,g) \in \mathcal{Q}_{ad} \times H^{-1/2}_{-c_f}(\partial\Omega) \times H^{1/2}_\diamond(\partial\Omega)$ be arbitrary. Assume that $\mathcal{N}_{ q}j, \mathcal{D}_{ q}g \in H^2(\Omega)$. Then
\begin{equation}\label{eq:r}
\widehat{\varrho}^h_{q} \left(j,g\right)\le C_r \big( h|\log h|\big)^r \mbox{ with }
\begin{cases}
r<1/2&\mbox{ if }d=2 \mbox{~and~}\\
r=1/3&\mbox{ if }d=3.
\end{cases}
\end{equation}
\end{lemma}

\begin{proof}
Due to Lemma \ref{moli.data}, since $\mathcal{N}_{ q}j \in H^2(\Omega)$, we get that
\begin{align}\label{7-4-16ct7}
\left\| \mathcal{N}_{q}j - \Pi_\diamond^h \mathcal{N}_{q}j \right\|_{H^1(\Omega)} \le Ch.
\end{align}
Furthermore,  it follows from Lemma \ref{23-6-16ct5} that
\begin{align}\label{7-4-16ct8}
\big\| q - \widehat{q}^h\big\|_{L^p(\Omega)}
=\left(\int_\Omega |q-\widehat{q}^h|\, |q-\widehat{q}^h|^{p-1}\right)^{1/p}
\le \left((2\overline{q})^{p-1}\, Ch|\log h| \right)^{1/p} \le C\big(h|\log h|\big)^{1/p}
\end{align}
for $p\in[1,\infty)$.
Like in \eqref{decomplemlim}, using (\ref{10/4:ct1}) and \eqref{ct9}, we infer that
\begin{align*}
\int_\Omega  \widehat{q}^h \nabla\mathcal{N}^h_{\widehat{q}^h}j \cdot \nabla \varphi^h
= \big\langle j,\gamma \varphi^h\big\rangle + \left( f,\varphi^h\right)
= \int_\Omega  q \nabla\mathcal{N}_{q}j \cdot \nabla \varphi^h
\end{align*}
for all $\varphi^h\in \mathcal{V}_{1,\diamond}^h$ and obtain that
\begin{align}\label{24-6-16ct1}
\int_\Omega  \widehat{q}^h \nabla \left( \mathcal{N}^h_{\widehat{q}^h}j - \Pi_\diamond^h \mathcal{N}_{q}j\right)  \cdot \nabla \varphi^h
= \int_\Omega  \big( q - \widehat{q}^h\big) \nabla\mathcal{N}_{q}j \cdot \nabla \varphi^h  + \int_\Omega  \widehat{q}^h \nabla \left( \mathcal{N}_{q}j - \Pi_\diamond^h \mathcal{N}_{q}j\right)  \cdot \nabla \varphi^h.
\end{align}
Since $H^2(\Omega)$ is embedded in $W^{1,s}(\Omega)$ with
\begin{equation}\label{eqs}
s\begin{cases}
<\infty&\mbox{ if }d=2\\
=6&\mbox{ if }d=3
\end{cases}
\end{equation}
(see, for example, \cite[Theorem 5.4]{adams}), it follows from Cauchy-Schwarz and H\"older's inequality that
\begin{align*}
\int_\Omega  \big( q - \widehat{q}^h\big) \nabla\mathcal{N}_{q}j \cdot \nabla \varphi^h
&\le \left( \int_\Omega  \big( q - \widehat{q}^h\big)^2 |\nabla\mathcal{N}_{q}j|^2 \right)^{1/2} \left( \int_\Omega  \big|\nabla \varphi^h\big|^2 \right)^{1/2}\\
&\le \big\| q - \widehat{q}^h\big\|_{L^{2s/(s-2)}(\Omega)}
\big\|\nabla\mathcal{N}_{q}j\big\|_{L^s(\Omega)} \big\| \varphi^h\big\|_{H^1(\Omega)}\\
&\le C \big\| q - \widehat{q}^h\big\|_{L^{2s/(s-2)}(\Omega)} \big\| \varphi^h\big\|_{H^1(\Omega)}.
\end{align*}
Then taking $\varphi^h = \mathcal{N}^h_{\widehat{q}^h}j - \Pi_\diamond^h \mathcal{N}_{q}j \in \mathcal{V}_{1,\diamond}^h$ and using \eqref{coercivity}, we infer from \eqref{24-6-16ct1} that
\begin{align*}
\left\| \mathcal{N}^h_{\widehat{q}^h}j - \Pi_\diamond^h \mathcal{N}_{q}j\right\|_{H^1(\Omega)}
&\le C \left( \big\| q - \widehat{q}^h\big\|_{L^{2s/(s-2)}(\Omega)} + \left\| \mathcal{N}_{q}j - \Pi_\diamond^h \mathcal{N}_{q}j \right\|_{H^1(\Omega)} \right)\\
&\le C\big(h|\log h|\big)^{(s-2)/(2s)} + Ch
\ \le C\big(h|\log h|\big)^{(s-2)/(2s)},
\end{align*}
by \eqref{7-4-16ct7}--\eqref{7-4-16ct8}. Thus, applying the triangle inequality and \eqref{7-4-16ct7} again, we infer that
\begin{align*}
\left\| \mathcal{N}_{q}j - \mathcal{N}^h_{\widehat{q}^h}j \right\|_{H^1(\Omega)}
&\le \left\| \mathcal{N}_{q}j - \Pi_\diamond^h \mathcal{N}_{q}j \right\|_{H^1(\Omega)} + \left\| \Pi_\diamond^h \mathcal{N}_{q}j - \mathcal{N}^h_{\widehat{q}^h}j\right\|_{H^1(\Omega)} \le C\big(h|\log h|\big)^{(s-2)/(2s)}.
\end{align*}
Similarly, we also get
$\left\| \mathcal{D}_{q} g - \mathcal{D}^h_{\widehat{q}^h} g \right\|_{H^1(\Omega)}  \le C\big(h|\log h|\big)^{(s-2)/(2s)}$ and so that
$$\widehat{\varrho}^h_{q} \left(j,g\right)\le C\big(h|\log h|\big)^{(s-2)/(2s)}$$ for $s$ as in \eqref{eqs}, which yields the assertion.
\end{proof}

With an appropriate a priori choice of the regularization parameter we get convergence under conditions similar to those stated, e.g., in \cite{NeubauerScherzer89} in the Hilbert space setting.

\begin{theorem}\label{convergence1}
Let $\left(h_n\right)_n$,  $\left(\delta_n\right)_n$
and $\left(\rho_n\right)_n$ be any positive sequences such that
\begin{align}\label{29-6-16ct1}
\rho_n\rightarrow 0, ~\frac{\delta_n}{\sqrt{\rho_n}}
\rightarrow 0 \mbox{~and~} \frac{\widehat{\varrho}^{h_n}_{q} \left(j^\dag,g^\dag\right)}{\sqrt{\rho_n}}
\rightarrow 0 \mbox{~as~} n\to\infty,
\end{align}
where $q$ is any solution to $\mathcal{N}_q j^\dag = \mathcal{D}_q g^\dag$. Moreover, assume that $\big(j_{\delta_n}, g_{\delta_n}\big) $ is a sequence satisfying
$$\big\|j_{\delta_n} - j^\dag\big\|_{H^{-1/2}(\partial\Omega)} + \big\|g_{\delta_n} - g^\dag\big\|_{H^{1/2}(\partial\Omega)} \le \delta_n$$
and that $q_n := q_{\rho_n, \delta_n}^{h_n}$ is an arbitrary minimizer of $\left( \mathcal{P}_{\rho_n,\delta_n}^{h_n} \right)$ for each $n\in\N$.
Then a subsequence of $(q_n)$ which is not relabelled and a solution $q^\dag$ to
$\left( \mathcal{IP}\right)$ exist such that
\begin{equation}\label{eq:conv}
\limn \big\| q_n - q^\dag \big \|_{L^1(\Omega)} = 0 \mbox{~and~} \limn \int_\Omega |\nabla q_n | = \int_\Omega \big| \nabla q^\dag \big|.
\end{equation}
Furthermore, $\left( \mathcal{N}^{h_n}_{q_n}j_{\delta_n}\right) $ and $\left( \mathcal{D}^{h_n}_{q_n}g_{\delta_n}\right) $ converge to the unique weak solution $\Phi^\dag = \Phi^\dag(q^\dag,j^\dag,g^\dag)$ of the boundary value problem \eqref{17-5-16ct1}--\eqref{17-5-16ct3} in the $H^1(\Omega)$-norm.
If $q^\dag$ is unique, then convergence \eqref{eq:conv} holds for the whole sequence.
\end{theorem}

Uniform $L^\infty$ boundedness of $(q_n)$ together with interpolation implies that convergence actually takes place in any $L^p$ space with $p\in[1,\infty]$.

\begin{remark}
In case $\mathcal{N}_{ q}j^\dag, \mathcal{D}_{ q}g^\dag \in H^2(\Omega)$ Lemma \ref{7-4-16ct3} shows that $\widehat{\varrho}^h_{q} \left(j^\dag,g^\dag\right)\le C\big(h|\log h|\big)^r$ with $r$ as in \eqref{eq:r}. Therefore, in view of \eqref{29-6-16ct1}, convergence is obtained if the sequence $(\rho_n)$ is chosen such that
\begin{align*}
\rho_n\rightarrow 0, ~\frac{\delta_n}{\sqrt{\rho_n}}
\rightarrow 0 \mbox{~and~} \frac{\big(h_n|\log h_n|\big)^r}{\sqrt{\rho_n}}
\rightarrow 0 \mbox{~as~} n\to\infty.
\end{align*}
By regularity theory for elliptic boundary value problems (see, for example, \cite{Grisvad,Troianiello}), if $j^\dag \in H^{1/2}(\Omega)$, $g^\dag \in H^{3/2}(\Omega)$, $q\in C^{0,1}(\Omega)$, $f\in L^2(\Omega)$ and either $\partial\Omega$ is $C^{1,1}$-smooth or the domain $\Omega$ is convex, then $\mathcal{N}_{ q}j^\dag, \mathcal{D}_{ q}g^\dag \in H^2(\Omega)$.
\end{remark}

\begin{proof}[Proof of Theorem \ref{convergence1}]
We have from the optimality of $q_n$ that
\begin{align}\label{eq:opt}
\mathcal{J}^{h_n}_{\delta_n} \left(q_n\right) + \rho_n \int_\Omega \sqrt{| \nabla q_n |^2 + \epsilon^{h_n}}&\le \mathcal{J}^{h_n}_{\delta_n} \big({\widehat{q}}^{h_n}\big) + \rho_n \int_\Omega \sqrt{\big| \nabla {\widehat{q}}^{h_n} \big|^2 + \epsilon^{h_n}},
\end{align}
where $ {\widehat{q}}^{h_n}$ is generated from $q$ according to Lemma \ref{23-6-16ct5}, and
\begin{align*}
&\mathcal{J}^{h_n}_{\delta_n} \big({\widehat{q}}^{h_n}\big)
= \int_\Omega {\widehat{q}}^{h_n} \nabla \left(\mathcal{N}^{h_n}_{{\widehat{q}}^{h_n}}j_{\delta_n}-
\mathcal{D}^{h_n}_{{\widehat{q}}^{h_n}}g_{\delta_n}\right) \cdot \nabla \left(\mathcal{N}^{h_n}_{{\widehat{q}}^{h_n}}j_{\delta_n}-
\mathcal{D}^{h_n}_{{\widehat{q}}^{h_n}}g_{\delta_n}\right)\\
&\quad\le \overline{q} \left\|\mathcal{N}^{h_n}_{{\widehat{q}}^{h_n}}j_{\delta_n}-
\mathcal{D}^{h_n}_{{\widehat{q}}^{h_n}}g_{\delta_n}\right\|^2_{H^1(\Omega)}\\
&\quad= \overline{q}\left\|\mathcal{N}^{h_n}_{{\widehat{q}}^{h_n}} j_{\delta_n} - \mathcal{N}^{h_n}_{{\widehat{q}}^{h_n}} j^\dag
+ \mathcal{N}^{h_n}_{{\widehat{q}}^{h_n}} j^\dag - \mathcal{D}^{h_n}_{{\widehat{q}}^{h_n}} g^\dag
- \mathcal{N}_q j^\dag + \mathcal{D}_q g^\dag
+ \mathcal{D}^{h_n}_{{\widehat{q}}^{h_n}} g^\dag - \mathcal{D}^{h_n}_{{\widehat{q}}^{h_n}}g_{\delta_n} \right\|^2_{H^1(\Omega)}\\
&\quad\le 4\overline{q}\Bigl( \left\|\mathcal{N}^{h_n}_{{\widehat{q}}^{h_n}} j_{\delta_n} - \mathcal{N}^{h_n}_{{\widehat{q}}^{h_n}} j^\dag \right\|^2_{H^1(\Omega)}
+ \left\| \mathcal{D}^{h_n}_{{\widehat{q}}^{h_n}} g^\dag - \mathcal{D}^{h_n}_{{\widehat{q}}^{h_n}}g_{\delta_n} \right\|^2_{H^1(\Omega)}\\
&~\qquad + \left\|\mathcal{N}^{h_n}_{{\widehat{q}}^{h_n}} j^\dag - \mathcal{N}_{q} j^\dag  \right\|^2_{H^1(\Omega)}
+  \left\| \mathcal{D}^{h_n}_{{\widehat{q}}^{h_n}} g^\dag  -\mathcal{D}_{q} g^\dag\right\|^2_{H^1(\Omega)}\Bigr)\\
&\quad\le C \left( \left\| j_{\delta_n} - j^\dag \right\|^2_{H^{-1/2}(\partial\Omega)} + \left\| g_{\delta_n} - g^\dag \right\|^2_{H^{-1/2}(\partial\Omega)}\right) + C \widehat{\varrho}^{h_n}_{q} \left(j^\dag,g^\dag\right)^2\\
&\quad\le C\left( \delta^2_n + \widehat{\varrho}^{h_n}_{q} \left(j^\dag,g^\dag\right)^2\right),
\end{align*}
where we have used Lemma \ref{6-5-16ct5} and the fact $\mathcal{N}_{q} j^\dag = \mathcal{D}_{q} g^\dag$.
Moreover, by Lemma \ref{23-6-16ct5}, we have that
\begin{align}\label{eq:limsup}
\limn \int_\Omega \sqrt{\big|\nabla {\widehat{q}}^{h_n} \big|^2 + \epsilon^{h_n}} = \limn \int_\Omega \big|\nabla {\widehat{q}}^{h_n} \big|
= \int_\Omega |\nabla q |.
\end{align}
We therefore conclude from \eqref{eq:opt} and \eqref{29-6-16ct1} that
\begin{align}\label{30-3-16ct2}
\limn \frac{\mathcal{J}^{h_n}_{\delta_n} \left( \widehat{q}^{h_n}\right)}{\rho_n} =0, \quad
\limn \mathcal{J}^{h_n}_{\delta_n} \left( q_n\right) =0
\end{align}
and
\begin{align}\label{30-3-16ct3}
\limsupn \int_\Omega | \nabla q_n |
\le \limsupn \int_\Omega \sqrt{\big| \nabla q_n \big|^2 + \epsilon^{h_n}} \le \limsupn \int_\Omega \sqrt{\big|\nabla {\widehat{q}}^{h_n} \big|^2 + \epsilon^{h_n}}
= \int_\Omega |\nabla q |.
\end{align}
Thus, $\left(q_n\right)$ is bounded in the $BV(\Omega)$-norm. A subsequence which is denoted by the same symbol and an element $q^\dag\in \mathcal{Q}_{ad}$ exist such that $\left(q_n\right)$ converges to $q^\dag$ in the
$L^1(\Omega)$-norm and
\begin{align}\label{30-3-16ct4}
\int_\Omega \big|\nabla q^\dag\big| \le \liminfn \int_\Omega |\nabla q_n |.
\end{align}
Using Lemma \ref{6-5-16ct5} again, we infer that
\begin{align*}
\Big\|\mathcal{N}^{h_n}_{q_n} j^\dag &- \mathcal{D}^{h_n}_{q_n} g^\dag  \Big\|^2_{H^1(\Omega)}  \\
&\le 3\left( \left\|\mathcal{N}^{h_n}_{q_n} j^\dag -  \mathcal{N}^{h_n}_{q_n}j_{\delta_n}  \right\|^2_{H^1(\Omega)} + \left\|\mathcal{D}^{h_n}_{q_n} g^\dag -  \mathcal{D}^{h_n}_{q_n}g_{\delta_n} \right\|^2_{H^1(\Omega)}  +\left\|\mathcal{N}^{h_n}_{q_n} j_{\delta_n} - \mathcal{D}^{h_n}_{q_n} g_{\delta_n}  \right\|^2_{H^1(\Omega)} \right) \notag\\
&\le C\delta^2_n +3\left\|\mathcal{N}^{h_n}_{q_n} j_{\delta_n} - \mathcal{D}^{h_n}_{q_n} g_{\delta_n}  \right\|^2_{H^1(\Omega)} \\
&\le C \left( \delta^2_n + \mathcal{J}^{h_n}_{\delta_n} \left( q_n\right)\right).
\end{align*}
Thus, using Lemma \ref{weakly conv.dis}, we obtain from the last inequality and \eqref{30-3-16ct2} that
\begin{align*}
\left\|\mathcal{N}_{q^\dag} j^\dag - \mathcal{D}_{q^\dag} g^\dag  \right\|^2_{H^1(\Omega)} = \limn \left\|\mathcal{N}^{h_n}_{q_n} j^\dag - \mathcal{D}^{h_n}_{q_n} g^\dag  \right\|^2_{H^1(\Omega)}  =0
\end{align*}
and so that
\begin{align}\label{eq:qdag}
\mathcal{N}_{q^\dag} j^\dag = \mathcal{D}_{q^\dag} g^\dag, \mbox{~i.e.,~} q^\dag \in \mathcal{I}_{\mathcal{Q}_{ad}} \left(j^\dag,g^\dag\right).
\end{align}
Furthermore, it follows from \eqref{30-3-16ct3}--\eqref{30-3-16ct4} that
$$\int_\Omega \big|\nabla q^\dag\big| \le \liminfn \int_\Omega |\nabla q_n| \le \limsupn \int_\Omega | \nabla q_n | \le \int_\Omega |\nabla q |,$$
for any solution $q$ to $\mathcal{N}_q j^\dag = \mathcal{D}_q g^\dag$,
hence, in view of \eqref{eq:qdag}, $q^\dag$ is a total variation minimizing solution of the identification problem, i.e., a solution to $\left( \mathcal{IP}\right)$.
Moreover, by setting $q=q^\dagger$, we get
$$\int_\Omega \big|\nabla q^\dag \big| = \limn \int_\Omega |\nabla q_n |. $$

Finally, Lemma \ref{weakly conv.dis} shows that the sequence $\left( \mathcal{N}^{h_n}_{q_n}j_{\delta_n}, \mathcal{D}^{h_n}_{q_n}g_{\delta_n}\right)$ converges in the $H^1(\Omega) \times H^1(\Omega)$-norm to $\left( \mathcal{N}_{q^\dag} j^\dag, \mathcal{D}_{q^\dag} g^\dag\right)$, where $\Phi^\dag := \mathcal{N}_{q^\dag} j^\dag = \mathcal{D}_{q^\dag} g^\dag$ is the unique weak solution of the elliptic system \eqref{17-5-16ct1}--\eqref{17-5-16ct3}, which finishes the proof.
\end{proof}

\section{Projected Armijo algorithm and numerical test}\label{11-7-16ct2}

In this section we present the {\it projected Armijo algorithm} (see \cite[Chapter 5]{kelley}) for numerically solving the minimization problem $\left(\mathcal{P}^h_{\rho,\delta}\right)$. We note that
many other efficient solution methods are available, see for example \cite{blank}.

\subsection{Projected Armijo algorithm}\label{PAA}

\subsubsection{Differentiability of the cost functional}

Similarly to Lemma \ref{q-diff} one also sees that the discrete Neumann and Dirichlet operators $\mathcal{N}^h$, $\mathcal{D}^h$ are Fr\'echet differentiable on the set $\mathcal{Q}$. For given $j_\delta\in H^{-1/2}(\partial\Omega)$ and each $q \in \mathcal{Q}$ the Fr\'echet derivative ${\mathcal{N}^h}'(q)\xi=:{\mathcal{N}^h_q}'j_\delta (\xi)$ in the direction $\xi \in L^\infty(\Omega)$ is an element of $\mathcal{V}^h_{1,\diamond}$ and satisfies the equation
\begin{align}
\int_\Omega q\nabla {\mathcal{N}^h_q}'j_\delta (\xi) \cdot \nabla \varphi^h
&= -\int_{\Omega} \xi\nabla \mathcal{N}^h_q j_\delta \cdot\nabla \varphi^h   \label{ct21***}
\end{align}
for all $\varphi^h \in \mathcal{V}^h_{1,\diamond}$.
Likewise, for fixed $g_\delta\in H^{1/2}(\partial\Omega)$ and each $q \in \mathcal{Q}$ the Fr\'echet derivative ${\mathcal{D}^h}'(q)\xi=:{\mathcal{D}^h_q}'g_\delta (\xi)$ in the direction $\xi \in L^\infty(\Omega)$  is an element of $\mathcal{V}_{1,0}^h$ and satisfies the equation
\begin{align}
\int_\Omega q\nabla {\mathcal{D}^h_q}'g_\delta (\xi) \cdot \nabla \psi^h
&= -\int_{\Omega} \xi\nabla \mathcal{D}^h_q g_\delta \cdot\nabla\psi^h \label{ct21****}
\end{align}
for all $\psi^h \in \mathcal{V}_{1,0}^h$.

The functional $\mathcal{J}_\delta^h$ is therefore Fr\'echet differentiable on the set $\mathcal{Q}$. For each $q\in \mathcal{Q}$ the action of the Fr\'echet derivative in the direction $ \xi \in L^\infty(\Omega)$ is given by
\begin{align*}
{\mathcal{J}_\delta^h}'(q) (\xi)
&= \int_\Omega \xi \nabla \left(\mathcal{N}^h_qj_\delta-
\mathcal{D}^h_qg_\delta\right) \cdot \nabla \left(\mathcal{N}^h_qj_\delta-
\mathcal{D}^h_qg_\delta\right) \\
&~\quad + 2\int_\Omega q \nabla \left({\mathcal{N}^h_q}'j_\delta (\xi)-
{\mathcal{D}^h_q}' g_\delta (\xi)\right) \cdot \nabla \left(\mathcal{N}^h_qj_\delta-
\mathcal{D}^h_qg_\delta\right)\\
&= \int_\Omega \xi \nabla \left(\mathcal{N}^h_qj_\delta-
\mathcal{D}^h_qg_\delta\right) \cdot \nabla \left(\mathcal{N}^h_qj_\delta-
\mathcal{D}^h_qg_\delta\right) \\
&~\quad + 2\int_\Omega q \nabla {\mathcal{N}^h_q}'j_\delta (\xi) \cdot \nabla \left(\mathcal{N}^h_qj_\delta-
\mathcal{D}^h_qg_\delta\right) -2\int_\Omega q \nabla \mathcal{N}^h_qj_\delta \cdot \nabla {\mathcal{D}^h_q}' g_\delta (\xi) \\
&~\quad + 2\int_\Omega q \nabla \mathcal{D}^h_qg_\delta \cdot \nabla {\mathcal{D}^h_q}' g_\delta (\xi).
\end{align*}
Since $\mathcal{N}^h_qj_\delta, \mathcal{D}^h_q g_\delta \in \mathcal{V}^h_{1,\diamond}$ and ${\mathcal{D}^h_q}' g_\delta (\xi) \in \mathcal{V}^h_{1,0} \subset \mathcal{V}^h_{1,\diamond}$, it follows from \eqref{ct21***}, \eqref{10/4:ct1} and \eqref{10/4:ct1*} that
\begin{align*}
\int_\Omega &q \nabla {\mathcal{N}^h_q}'j_\delta (\xi) \cdot \nabla \left(\mathcal{N}^h_qj_\delta-
\mathcal{D}^h_qg_\delta\right) -\int_\Omega q \nabla \mathcal{N}^h_qj_\delta \cdot \nabla {\mathcal{D}^h_q}' g_\delta (\xi) + \int_\Omega q \nabla \mathcal{D}^h_qg_\delta \cdot \nabla {\mathcal{D}^h_q}' g_\delta (\xi)\\
&= - \int_\Omega \xi \nabla \mathcal{N}^h_q j_\delta \cdot \nabla \left(\mathcal{N}^h_qj_\delta -\mathcal{D}^h_qg_\delta\right) -\left\langle  j_\delta, \gamma {\mathcal{D}^h_q}' g_\delta (\xi)\right\rangle - \left( f,  {\mathcal{D}^h_q}' g_\delta (\xi)\right) + \left( f,  {\mathcal{D}^h_q}' g_\delta (\xi)\right)\\
&= - \int_\Omega \xi \nabla \mathcal{N}^h_q j_\delta \cdot \nabla \left(\mathcal{N}^h_qj_\delta -\mathcal{D}^h_qg_\delta\right)
\end{align*}
and so that
\begin{align*}
{\mathcal{J}_\delta^h}'(q) (\xi)
&= \int_\Omega \xi \nabla \left(\mathcal{N}^h_qj_\delta-
\mathcal{D}^h_qg_\delta\right) \cdot \nabla \left(\mathcal{N}^h_qj_\delta-
\mathcal{D}^h_qg_\delta\right) -2\int_\Omega \xi \nabla \mathcal{N}^h_q j_\delta \cdot \nabla \left(\mathcal{N}^h_qj_\delta -\mathcal{D}^h_qg_\delta\right)\\
&=\int_\Omega \xi \left( \nabla\mathcal{D}^h_qg_\delta \cdot \nabla\mathcal{D}^h_qg_\delta - \nabla\mathcal{N}^h_qj_\delta \cdot \nabla\mathcal{N}^h_qj_\delta\right).
\end{align*}
Therefore, the derivative of the cost functional $\Upsilon^h_{\rho,\delta}$ of $\left(\mathcal{P}^h_{\rho,\delta}\right)$ at $q\in \mathcal{Q}^h_{ad}$ in the direction $ \xi \in \mathcal{V}^h_1$ is given by
\begin{align}\label{1-8-16ct1}
{\Upsilon^h_{\rho,\delta}}'(q) (\xi) = \int_\Omega \xi \left( \nabla\mathcal{D}^h_qg_\delta \cdot \nabla\mathcal{D}^h_qg_\delta - \nabla\mathcal{N}^h_qj_\delta \cdot \nabla\mathcal{N}^h_qj_\delta\right) + \rho  \int_\Omega \frac{\nabla q \cdot \nabla \xi}{\sqrt{|\nabla q|^2 +\epsilon^h}}.
\end{align}
Let $\{N_j~:~j=1,\ldots,M^h\}$ be the set of nodes of the triangulation $\mathcal{T}^h$, then $\mathcal{V}^h_1$ is a finite dimensional vector space with dimension $M^h$. Let $\{\phi_1, \ldots,\phi_{M^h}\}$ be the basis of $\mathcal{V}^h_1$ consisting hat functions, i.e. $\phi_i(N_j) =\delta_{ij}$ for all $1\le i,j\le M^h$, where $\delta_{ij}$ is the Kronecker symbol. Each functional $u\in\mathcal{V}^h_1$ then can be identified with a vector $(u_1,\ldots,u_{M^h})\in\mathbb{R}^{M^h}$ consisting of the nodal values of $u$, i.e.
$$u=\sum_{j=1}^{M^h} u_j \phi_j \mbox{~with~} u_j=u(N_j).$$
In $\mathcal{V}^h_1$ we use the Euclidean inner product $\langle \cdot,\cdot\rangle_E$. For each $u=(u_1,\ldots,u_{M^h})$ and $v=(v_1,\ldots,v_{M^h})$, we have
$\langle u,v\rangle_E =\sum_{j=1}^{M^h} u_jv_j$. Let us denote the gradient of $\Upsilon^h_{\rho,\delta}$ at $q\in \mathcal{Q}^h_{ad}$ by $\nabla \Upsilon^h_{\rho,\delta}(q) = (\Upsilon_1,\ldots,\Upsilon_{M^h})$. We then have from \eqref{1-8-16ct1} with $ \xi=(\xi_1,\ldots,\xi_{M^h}) \in \mathcal{V}^h_1$ that
$$\sum_{j=1}^{M^h} \xi_j \int_\Omega \left( \phi_j \left( \nabla\mathcal{D}^h_{q}g_\delta \cdot \nabla\mathcal{D}^h_{q}g_\delta - \nabla\mathcal{N}^h_{q}j_\delta \cdot \nabla\mathcal{N}^h_{q}j_\delta\right) + \frac{\rho \nabla q \cdot \nabla \phi_j}{\sqrt{|\nabla q|^2 +\epsilon^h}}\right) =\sum_{j=1}^{M^h} \xi_j \Upsilon_j$$
which yields
\begin{align}\label{1-8-16ct2}
\Upsilon_j = \int_\Omega \phi_j \left( \nabla\mathcal{D}^h_{q}g_\delta \cdot \nabla\mathcal{D}^h_{q}g_\delta - \nabla\mathcal{N}^h_{q}j_\delta \cdot \nabla\mathcal{N}^h_{q}j_\delta\right) + \rho\int_\Omega \frac{\nabla q \cdot \nabla \phi_j}{\sqrt{|\nabla q|^2 +\epsilon^h}}
\end{align}
for all $j=1,\ldots,M^h$.

\subsubsection{Algorithm}

The projected Armijo algorithm is then read as: given a step size control $\beta \in (0,1)$, an initial approximation $q^h_0 \in \mathcal{Q}^h_{ad}$, a smoothing parameter $\epsilon^h$, number of iteration $N$ and setting $k=0$.
\begin{enumerate}
\item Compute $\mathcal{N}^h_{q^h_k} j_\delta$ and $\mathcal{D}^h_{q^h_k} g_\delta$ from the variational equations
\begin{align}\label{eqNh}
\int_\Omega q^h_k \nabla \mathcal{N}^h_{q^h_k} j_\delta \cdot \nabla \varphi^h = \big\langle j_\delta,\gamma\varphi^h\big\rangle + \left( f,\varphi^h\right) \mbox{~for all~} \varphi^h\in
\mathcal{V}_{1,\diamond}^h
\end{align}
and
\begin{align}\label{eqDh}
\int_\Omega q^h_k \nabla \mathcal{D}^h_{q^h_k} g_\delta \cdot \nabla \psi^h = \left( f,\psi^h\right) \mbox{~for all~} \psi^h\in
\mathcal{V}_{1,0}^h,
\end{align}
respectively, as well as $\Upsilon_{\rho,\delta}^h (q^h_k)$ according to \eqref{27-3-16ct4}, \eqref{29/6:ct9}.
\item Compute the gradient $\nabla \Upsilon^h_{\rho,\delta}(q^h_k)$ with the $j^{\mbox {th}}$-component given by
$$\Upsilon_j = \int_\Omega \phi_j \left( \nabla\mathcal{D}^h_{q^h_k}g_\delta \cdot \nabla\mathcal{D}^h_{q^h_k}g_\delta - \nabla\mathcal{N}^h_{q^h_k}j_\delta \cdot \nabla\mathcal{N}^h_{q^h_k}j_\delta\right) + \rho  \int_\Omega \frac{\nabla q^h_k  \cdot \nabla \phi_j}{\sqrt{|\nabla q^h_k |^2 +\epsilon^h}},$$
due to \eqref{1-8-16ct2}.
\item Set $G^h_k :=\sum_{j=1}^{M^h} \Upsilon_j\phi_j$.
\begin{enumerate}
\item  Compute
$$\tilde{q}^h_k := \max\left( \underline{q}, \min\left(q^h_k - \beta G^h_k, \overline{q} \right) \right) ,$$
$\mathcal{N}^h_{\tilde{q}^h_k} j_\delta$, $\mathcal{D}^h_{\tilde{q}^h_k} g_\delta$, according to \eqref{eqNh}, \eqref{eqDh},
$\Upsilon_{\rho,\delta}^h (\tilde{q}^h_k)$,  according to \eqref{27-3-16ct4}, \eqref{29/6:ct9}, and \\
$L:= \Upsilon_{\rho,\delta}^h (\tilde{q}^h_k) - \Upsilon_{\rho,\delta}^h (q^h_k) + \tau\beta \big\| \tilde{q}^h_k - q^h_k \big\|^2_{L^2(\Omega)}
$ with $\tau=10^{-4}$.
\item If $L\le 0$

\qquad  go to the next step (c) below

else

\qquad  set $\beta := \frac{\beta}{2}$ and then go back (a)
\item Update $q^h_k = \tilde{q}^h_k$, set $k=k+1$.
\end{enumerate}
\item  Compute
\begin{align}\label{7-7-16ct1}
\mbox{Tolerance}:= \big\| \nabla \Upsilon^h_{\rho,\delta}(q^h_k) \big\|_{L^2(\Omega)} -\tau_1 -\tau_2\big\| \nabla \Upsilon^h_{\rho,\delta}(q^h_0) \big\|_{L^2(\Omega)}
\end{align}
with $\tau_1 := 10^{-3}h^{1/2}$ and $\tau_2 := 10^{-2}h^{1/2}$. If $\mbox{Tolerance} \le 0$ or $k>N$, then stop; otherwise go back Step 1.
\end{enumerate}

\subsection{Numerical tests}

We now illustrate the theoretical result with numerical examples.
For this purpose we consider the the boundary value problem
\begin{align}
-\nabla \cdot \big(q^\dag \nabla \Phi \big) &= f \mbox{~in~} \Omega,  \label{17-5-16ct1*}\\
q^\dag \nabla \Phi \cdot \vec{n} &= j^\dag \mbox{~on~} \partial\Omega \mbox{~and~} \label{17-5-16ct2*}\\
\Phi &= g^\dag \mbox{~on~} \partial\Omega \label{17-5-16ct3*}
\end{align}
with $\Omega = \{ x = (x_1,x_2) \in \R^2 ~|~ -1 < x_1, x_2 < 1\}$. The special constants $\underline{q}$ and $\overline{q}$ in the definition of the set $\mathcal{Q}$ according to \eqref{22-3-15ct3} are respectively chosen as 0.05 and 10.

We assume that the known source $f$ is discontinuous and given by
$$f= \frac{3}{2}\chi_D - \frac{1}{2} \chi_{\Omega\setminus D},$$
where $\chi_D$ is the characteristic function of $D:= \left\{ (x_1, x_2) \in \Omega ~\big|~ |x_1| \le 1/2 \mbox{~and~} |x_2| \le 1/2 \right\}$. Note that $(f,1)=0$, so that $c_f =0$. The sought conductivity $q^\dag$ in the equation \eqref{17-5-16ct1*}--\eqref{17-5-16ct2*} is assumed to be discontinuous and given by
$$q^\dag = 3 \chi_{\Omega_{1}} + 2 \chi_{\Omega_{2}} + \chi_{\Omega \setminus (\Omega_1 \cup \Omega_2)},$$
where
\begin{align*}
&\Omega_1 := \left\{ (x_1, x_2) \in \Omega ~\big|~ 9\Big(x_1+1/2\Big)^2 + 16\Big(x_2-1/2\Big)^2 \le 1\right\} \mbox{~and~}\\
&\Omega_2 := \left\{ (x_1, x_2) \in \Omega ~\big|~ \enskip\Big(x_1-1/2\Big)^2 + \Big(x_2+1/2\Big)^2 \le 1/16\right\}.
\end{align*}
For the discretization we divide the interval $(-1,1)$ into $\ell$ equal segments and so that the domain $\Omega = (-1,1)^2$ is divided into $2\ell^2$ triangles, where the diameter of each triangle is $h_{\ell} = \frac{\sqrt{8}}{\ell}$. In the minimization problem $\left(\mathcal{P}_{\rho,\delta}^{h} \right)$ we take $h=h_\ell$ and $\rho = \rho_\ell = 0.01 \sqrt{h_\ell}$. We use the projected Armijo algorithm which is described in Subsection \ref{PAA} for computing the numerical solution of the problem $\left(\mathcal{P}_{\rho_\ell,\delta_\ell}^{h_\ell} \right)$. The step size control is chosen with $\beta=0.75$ while the smoothing parameter $\epsilon^{h_\ell} = \rho_\ell$. The initial approximation is the constant function defined by $q^{h_\ell}_0 = 1.5$.

\begin{example}\label{30-6-17ct1} In this example the Neumann boundary condition $j^\dag \in H^{-1/2}_{-c_f}(\partial\Omega)$ in the equation \eqref{17-5-16ct2*} is chosen to be the piecewise constant function defined by
\begin{equation}\label{20-6-17ct5}
\begin{aligned}
j^\dag 
&= \chi_{(0,1]\times\{-1\}} - \chi_{[-1,0]\times\{1\}} + 2\chi_{(0,1]\times\{1\}} -2\chi_{[-1,0]\times\{-1\}}\\
&~\quad +3\chi_{\{-1\}\times(-1,0]} - 3\chi_{\{1\}\times(0,1)} + 4\chi_{\{1\}\times(-1,0]}  - 4\chi_{\{-1\}\times(0,1)} .
\end{aligned}
\end{equation}
so that
$\left\langle j^\dag, 1\right\rangle =0$. 
The Dirichlet boundary condition $g^\dag \in H^{1/2}_\diamond(\partial\Omega)$ in the equation \eqref{17-5-16ct3*} is then defined as
$g^\dag = \gamma \mathcal{N}_{q^\dag} j^\dag,$
where $\mathcal{N}_{q^\dag} j^\dag$ is the unique weak solution to the Neumann problem \eqref{17-5-16ct1*}--\eqref{17-5-16ct2*}. For the numerical solution of the pure Neumann problem \eqref{17-5-16ct1*}--\eqref{17-5-16ct2*} we use the penalty technique, see e.g. \cite{bochev,gockenbach} for more details. Furthermore, to avoid a so-called inverse crime, we generate the data on a finer grid than those used in the computations. To do so, we first solve the Neumann problem \eqref{17-5-16ct1*}--\eqref{17-5-16ct2*} on the very fine grid $\ell = 128$, and then handle $(j^\dag,g^\dag)$ on this grid for our computational process below.

We assume that noisy observations are available in the form
\begin{align}\label{3-7-17ct1}
\left( j_{\delta_{\ell}}, g_{\delta_{\ell}} \right) = \left( j^\dag+ \theta_\ell\cdot R_{j^\dag}, g^\dag+ \theta_\ell\cdot R_{g^\dag}\right) \quad \mbox{for some} \quad \theta_\ell>0 \quad \mbox{depending on} \quad \ell,
\end{align}
where $R_{j^\dag}$ and $R_{g^\dag}$ are $\partial M^{h_\ell}\times 1$-matrices of random numbers on the interval $(-1,1)$ which are generated by the MATLAB function ``rand'' and $\partial M^{h_\ell}$ is the number of boundary nodes of the triangulation $\mathcal{T}^{h_\ell}$. The measurement error is then computed as $\delta_\ell = \big\|j_{\delta_\ell} -j^\dag\big\|_{L^2(\partial\Omega)} + \big\|g_{\delta_\ell} -g^\dag\big\|_{L^2(\partial\Omega)}.$ 
To satisfy the condition $\delta_\ell\cdot\rho^{-1/2}_\ell\to 0$ as $\ell\to\infty$ in Theorem \ref{convergence1} we below take $\theta_\ell = h_\ell \sqrt{\rho_\ell}$.
In doing so, we reversely mimic the situation of a given sequence of noise levels $\delta_\ell$ tending to zero and of choosing the discretization level as well as the regularization parameter in dependence of the noise level.

Our computational process will be started with the coarsest level $\ell=4$. In each iteration $k$ we compute Tolerance
defined by \eqref{7-7-16ct1}. Then the iteration is stopped if
$\mbox{Tolerance} \le 0$ or the number of iterations reaches the maximum iteration count of 1000.
After obtaining the numerical solution of the first iteration process with respect to the coarsest level $\ell=4$, we use its interpolation on the next finer mesh $\ell=8$ as an initial approximation $q^{h_\ell}_0$ for the algorithm on this finer mesh, and so on for $\ell=16,32,64$.

Let $q_\ell$ denote the conductivity obtained at {\it the final iterate} of the algorithm corresponding to the refinement level $\ell$. Furthermore, let $\mathcal{N}^{h_\ell}_{q_\ell} j_{\delta_\ell}$ and $\mathcal{D}^{h_\ell}_{q_\ell} g_{\delta_\ell}$ denote {\it the computed numerical solution} to the Neumann problem
\begin{align*}
-\nabla \cdot (q_\ell \nabla u) = f \mbox{~in~} \Omega \mbox{~and~}
q_\ell \nabla u\cdot \vec{n} = j_{\delta_\ell} \mbox{~on~} \partial\Omega
\end{align*}
and the Dirichlet problem
\begin{align*}
-\nabla \cdot (q_\ell \nabla v) = f \mbox{~in~} \Omega \mbox{~and~}
v = g_\ell \mbox{~on~} \partial\Omega,
\end{align*}
respectively. The notations $\mathcal{N}^{h_\ell}_{q^\dag}  j^\dag$ and $\mathcal{D}^{h_\ell}_{q^\dag} g^\dag$ of {\it the exact numerical solutions} are to be understood similarly. We use the following abbreviations for the errors
\begin{align*}
L^2_q = \big\|q_\ell -  q^\dag\big\|_{L^2(\Omega)}, \quad L^2_\mathcal{N} = \big\|\mathcal{N}^{h_\ell}_{q_\ell} j_{\delta_\ell} - \mathcal{N}^{h_\ell}_{q^\dag}  j^\dag\big\|_{L^2(\Omega)} \quad \mbox{and} \quad L^2_\mathcal{D} = \big\|\mathcal{D}^{h_\ell}_{q_\ell} g_{\delta_\ell} - \mathcal{D}^{h_\ell}_{q^\dag}  g^\dag\big\|_{L^2(\Omega)}.
\end{align*}
The numerical results are summarized in Table \ref{b1} and Table \ref{b3}, where we present the refinement level $\ell$, the mesh size $h_\ell$ of the triangulation, the regularization parameter $\rho_\ell$, the measurement noise $\delta_\ell$, the number of iterations, the value of Tolerance, the errors $L^2_q$, $L^2_{\mathcal{N}}$, $L^2_{\mathcal{D}}$, and their experimental order of convergence (EOC) defined by
$$\mbox{EOC}_\Xi := \dfrac{\ln \Xi(h_1) - \ln \Xi(h_2)}{\ln h_1 - \ln h_2}$$
with $\Xi(h)$ being an error functional with respect to the mesh size $h$.
The convergence history given in Table \ref{b1} and Table \ref{b3} shows that the projected Armijo algorithm performs well for our identification problem.

All figures presented hereafter correspond to the finest level $\ell = 64$. Figure  \ref{h4} from left to right shows the interpolation $I_1^{h_{\ell}} q^\dag$, the numerical solution $q_\ell$  computed by the algorithm at the 953$^{\mbox{\tiny th}}$ iteration, and the differences $\mathcal{N}^{h_\ell}_{q_\ell} j_{\delta_\ell} - \mathcal{N}^{h_\ell}_{q^\dag}  j^\dag$ and $\mathcal{D}^{h_\ell}_{q_\ell} g_{\delta_\ell} - \mathcal{D}^{h_\ell}_{q^\dag} g^\dag$.

\begin{table}[H]
\begin{center}
\begin{tabular}{|c|l|l|l|r|l|}
\hline \multicolumn{6}{|c|}{ {\bf Convergence history} }\\
\hline
$\ell$ &\scriptsize $h_\ell$ &\scriptsize $\rho_\ell$ &\scriptsize $\delta_\ell$ &\scriptsize {\bf Iterate} &\scriptsize {\bf Tolerance} \\
\hline
4   &0.7071 & 8.4091e-3& 0.1733& 1000&  ~0.2459\\
\hline
8  &0.3536 & 5.9460e-3& 8.4273e-2& 1000&  ~3.2771e-2\\
\hline
16  &0.1766 & 4.2045e-3& 3.4320e-2& 1000&   ~5.8479e-3\\
\hline
32  &8.8388e-2 & 2.9730e-3& 1.5877e-2& 1000&  ~9.4359e-5\\
\hline
64  &4.4194e-2 & 2.1022e-3& 6.7743e-3& 953& -7.6116e-5\\
\hline
\end{tabular}
\caption{Refinement level $\ell$, mesh size $h_\ell$ of the triangulation, regularization parameter $\rho_\ell$,  measurement noise $\delta_\ell$, number of iterates and value of Tolerance.}
\label{b1}
\end{center}
\end{table}

\begin{table}[H]
\begin{center}
\begin{tabular}{|c|l|l|l|l|l|l|}
 \hline 
 \multicolumn{7}{|c|}{ {\bf Convergence history and EOC} }\\
 \hline
$\ell$ &\scriptsize $L^2_q$ &\scriptsize  $L^2_{\mathcal{N}}$ &\scriptsize $L^2_{\mathcal{D}}$ &\scriptsize EOC$_{L^2_q}$ &\scriptsize  EOC$_{L^2_{\mathcal{N}}}$ &\scriptsize EOC$_{L^2_{\mathcal{D}}}$\\
\hline
4&    0.7906 & 0.3016 & 0.1371 &--- &--- &---\\
\hline
8&   0.4768 & 0.1546 & 6.2771e-2 & 0.7296 & 0.9637  & 1.1271\\
\hline
16&   0.2306 & 6.9702e-2 & 2.1228e-2 & 1.0480 & 1.1497  & 1.5641\\
\hline
32&   0.1271 & 3.0668e-2 & 9.9234e-3 &0.8594  & 1.1845  & 1.0971\\
\hline
64&   6.7791e-2&  1.2116e-2 & 5.1055e-3 & 0.9068
 & 1.3398 & 0.9588\\
\hline
\multicolumn{4}{|c|}{\bf Mean of EOC} & {\bf 0.8859
} & {\bf 1.1594}  & {\bf 1.1868}\\
\hline
\end{tabular}
\caption{Errors $L^2_q$, $L^2_{\mathcal{N}}$, $L^2_{\mathcal{D}}$, and their EOC between finest and coarsest level.}
\label{b3}
\end{center}
\end{table}

\begin{figure}[H]
\begin{center}
\includegraphics[scale=0.20]{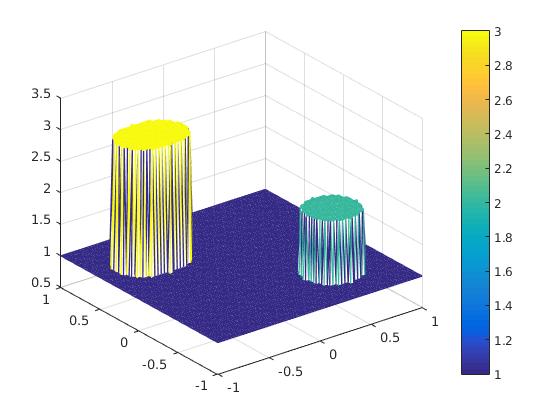}
\includegraphics[scale=0.20]{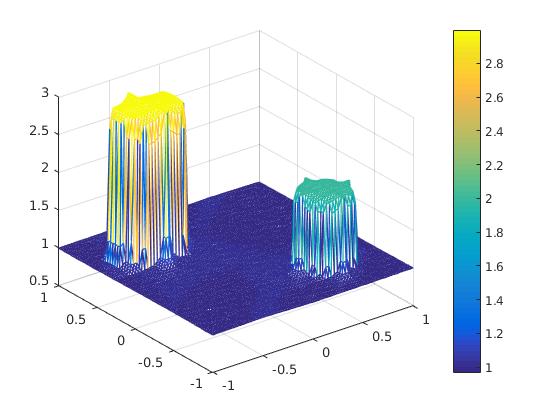} 
\includegraphics[scale=0.20]{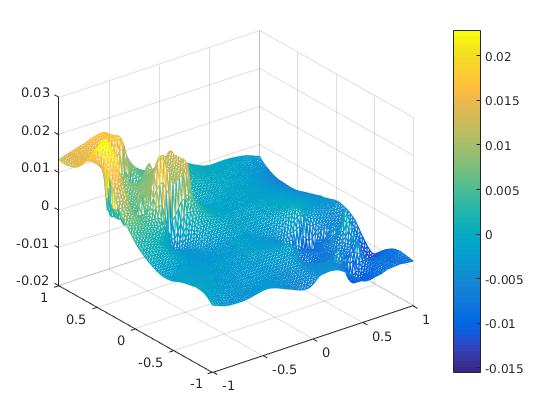}
\includegraphics[scale=0.20]{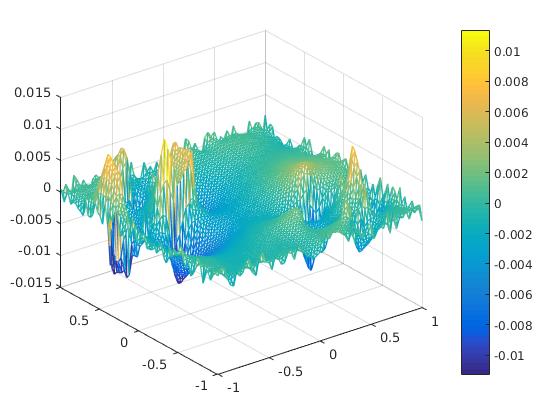}
\end{center}
\caption{Interpolation $I_1^{h_{\ell}} q^\dag$, computed numerical solution $q_\ell$ of the algorithm at the 953$^{\mbox{\tiny th}}$ iteration, and the differences $\mathcal{N}^{h_\ell}_{q_\ell} j_{\delta_\ell} - \mathcal{N}^{h_\ell}_{q^\dag}  j^\dag$ and $\mathcal{D}^{h_\ell}_{q_\ell} g_{\delta_\ell} - \mathcal{D}^{h_\ell}_{q^\dag} g^\dag$, for $\ell=64$, $\delta_\ell=6.7743e-3$.}
\label{h4}
\end{figure}
We observe a decrease of all errors as the noise level gets smaller, as expected from our convergence result, however, with respect to different norms. In particular, in our computations we use an $L^2$ noise level, as realistic in applications.
\end{example}

\begin{example}
In this example we consider noisy observations in the form
$$\left( j_{\delta_{\ell}}, g_{\delta_{\ell}} \right) = \left( j^\dag+ \theta\cdot R_{j^\dag}, g^\dag+ \theta\cdot R_{g^\dag}\right),$$
where $j^\dag$ is defined by \eqref{20-6-17ct5}. This is different from \eqref{3-7-17ct1},
since here $\theta>0$ is independent of $\ell$.

Using the computational process which was described as in Example \ref{30-6-17ct1} starting with $\ell=4$, in Table \ref{b4} we perform the numerical results for the finest grid $\ell=64$ and with different values of $\theta$.

\begin{table}[H]
\begin{center}
\begin{tabular}{|l|l|l|l|l|l|l|}
\hline \multicolumn{7}{|c|}{ {\bf Numerical results for the finest grid $\ell=64$} }\\
\hline
$\theta$ &\scriptsize $\delta_\ell$  &\scriptsize {\bf Iterate} &\scriptsize {\bf Tolerance} &\scriptsize $L^2_q$ &\scriptsize $L^2_\mathcal{N}$ &\scriptsize $L^2_\mathcal{D}$\\
\hline
0.005   &0.0167 & 991& -9.7239e-5& 7.7012e-2&  1.3085e-2& 7.5847e-3\\
\hline
0.01  &0.0316 & 1000& ~3.4160e-4& 9.7971e-2&  1.7896e-2& 9.3278e-3\\
\hline
0.05  &0.1567 & 1000& ~7.2599e-3& 0.2467&   0.1071 &3.8046e-2 \\
\hline
0.1  &0.3308 & 1000& ~2.2546e-2& 0.4059&  0.2067 &0.1077 \\
\hline
\end{tabular}
\caption{Numerical results for the finest grid $\ell=64$ and with different values of $\theta$.}
\label{b4}
\end{center}
\end{table}

In Figure \ref{h5} from left to right we show the computed numerical solution $q_\ell$ of the algorithm at the final iteration, and the differences $q_\ell-I^{h_\ell}_1q^\dag$, 
$\mathcal{N}^{h_\ell}_{q_\ell} j_{\delta_\ell} - \mathcal{N}^{h_\ell}_{q^\dag}  j^\dag$ and $\mathcal{D}^{h_\ell}_{q_\ell} g_{\delta_\ell} - \mathcal{D}^{h_\ell}_{q^\dag} g^\dag$ for $\ell=64$ and $\theta=0.005$, i.e., $\delta_\ell=0.0167$. Finally, Figure \ref{h6} performs the analog differences, but with $\theta=0.1$, i.e., $\delta_\ell=0.3308$.

\begin{figure}[H]
\begin{center}
\includegraphics[scale=0.19]{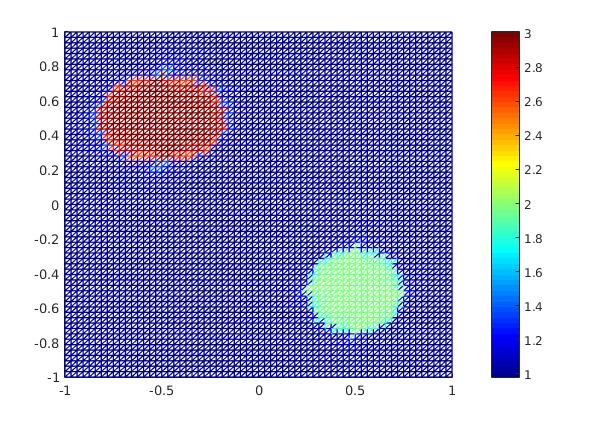}
\includegraphics[scale=0.19]{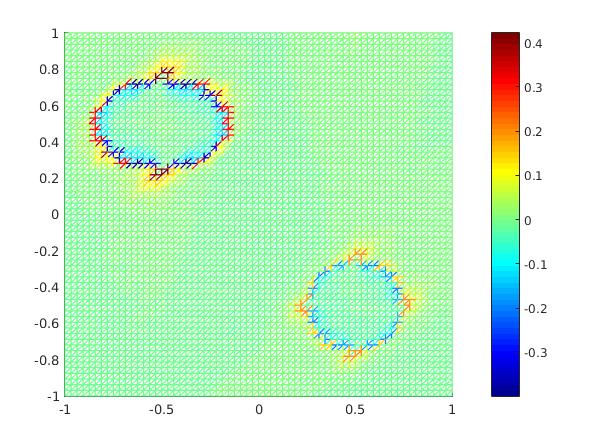} 
\includegraphics[scale=0.19]{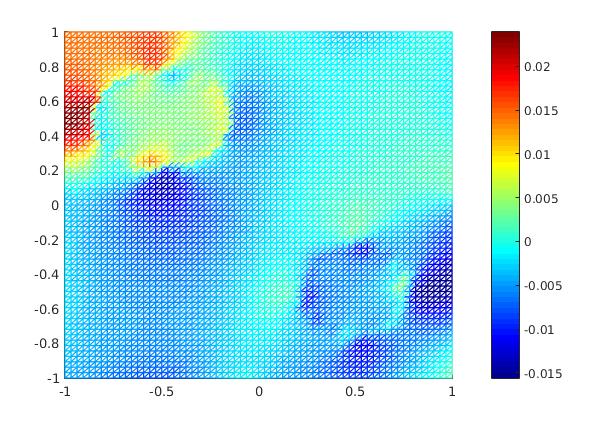}
\includegraphics[scale=0.19]{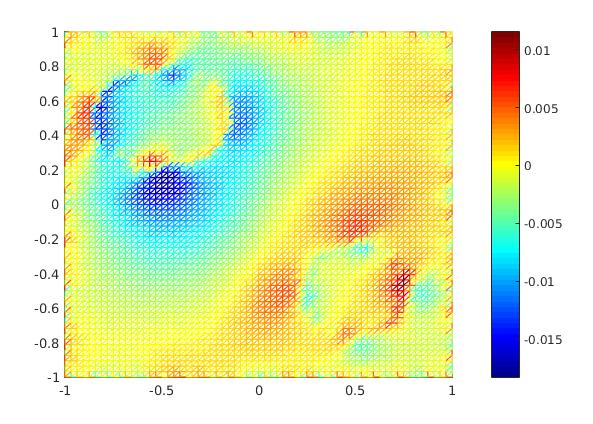}
\end{center}
\caption{Computed numerical solution $q_\ell$ of the algorithm at the 991$^{\mbox{\tiny th}}$ iteration, and the differences $q_\ell-I^{h_\ell}_1q^\dag$, 
$\mathcal{N}^{h_\ell}_{q_\ell} j_{\delta_\ell} - \mathcal{N}^{h_\ell}_{q^\dag}  j^\dag$ and $\mathcal{D}^{h_\ell}_{q_\ell} g_{\delta_\ell} - \mathcal{D}^{h_\ell}_{q^\dag} g^\dag$ for $\ell=64$ and $\theta=0.005$, i.e., $\delta_\ell=0.0167$.}
\label{h5}
\end{figure}

\begin{figure}[H]
\begin{center}
\includegraphics[scale=0.19]{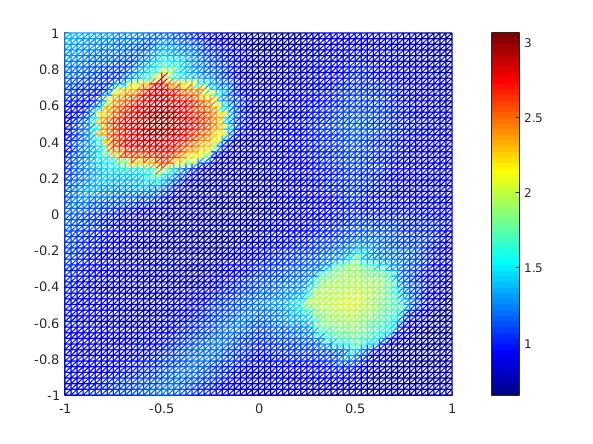}
\includegraphics[scale=0.19]{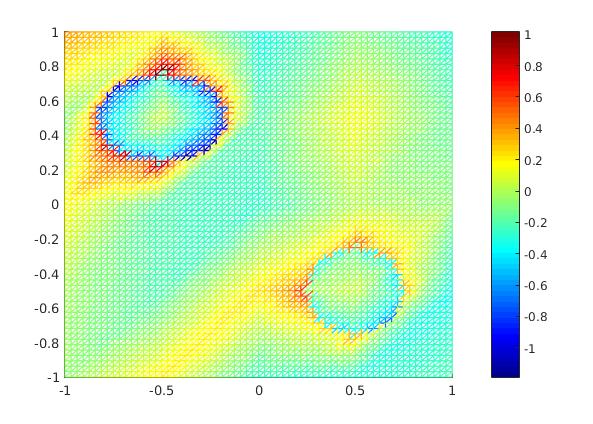} 
\includegraphics[scale=0.19]{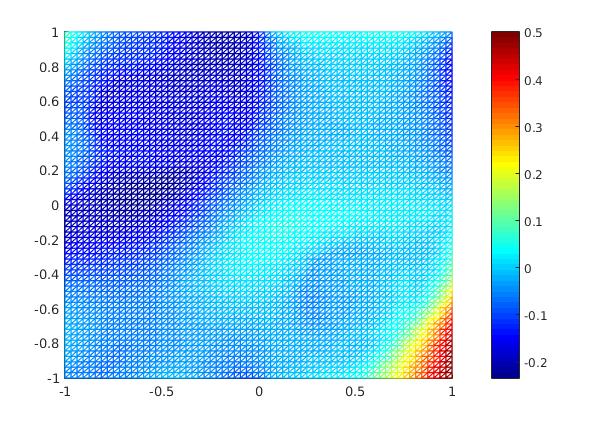}
\includegraphics[scale=0.19]{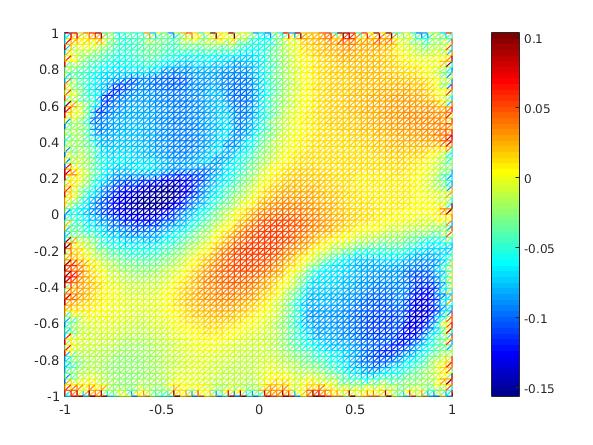}
\end{center}
\caption{Computed numerical solution $q_\ell$ of the algorithm at the 1000$^{\mbox{\tiny th}}$ iteration, and the differences $q_\ell-I^{h_\ell}_1q^\dag$, 
$\mathcal{N}^{h_\ell}_{q_\ell} j_{\delta_\ell} - \mathcal{N}^{h_\ell}_{q^\dag}  j^\dag$ and $\mathcal{D}^{h_\ell}_{q_\ell} g_{\delta_\ell} - \mathcal{D}^{h_\ell}_{q^\dag} g^\dag$ for $\ell=64$ and $\theta=0.1$, i.e., $\delta_\ell=0.3308$.}
\label{h6}
\end{figure}
\end{example}

\begin{example}\label{ex_mult}
In this example we assume that multiple measurements are available, say $\left(j_\delta^i,g_\delta^i \right)_{i=1,\ldots,I}$. Then, the cost functional $\Upsilon^h_{\rho,\delta}$ and the problem $\left(\mathcal{P}^h_{\rho,\delta}\right)$ can be rewritten as 
$$\min_{q\in\mathcal{Q}^h_{ad}} \bar{\Upsilon}^h_{\rho,\delta} (q) := \min_{q\in\mathcal{Q}^h_{ad}} \left( \underbrace{\frac{1}{I}\sum_{i=1}^I\int_\Omega q\nabla \left(\mathcal{N}^h_qj_\delta^i-
\mathcal{D}^h_qg_\delta^i\right) \cdot \nabla \left(\mathcal{N}^h_qj_\delta^i-
\mathcal{D}^h_qg_\delta^i\right)}_{:= \bar{\mathcal{J}}_\delta^h (q)} + \rho \int_\Omega \sqrt{\left|\nabla q\right|^2+\epsilon^h}\right),  \eqno\left(\bar{\mathcal{P}}^h_{\rho,\delta}\right)$$
which also attains a solution $\bar{q}^h_{\rho,\delta}$. The Neumann boundary condition in the equation \eqref{17-5-16ct2*} is chosen in the same form of \eqref{20-6-17ct5}, i.e.
\begin{equation}\label{20-6-17ct1}
\begin{aligned}
j^\dag_{(A,B,C,D)} 
&= A\cdot \chi_{(0,1]\times\{-1\}} - A\cdot\chi_{[-1,0]\times\{1\}} + B\cdot\chi_{(0,1]\times\{1\}} -B\cdot\chi_{[-1,0]\times\{-1\}} \\
&~\quad +C\cdot\chi_{\{-1\}\times(-1,0]} - C\cdot\chi_{\{1\}\times(0,1)}  + D\cdot\chi_{\{1\}\times(-1,0]} - D\cdot\chi_{\{-1\}\times(0,1)},
\end{aligned}
\end{equation}
that depends on the constants $A,B,C$ and $D$. Let $g^\dag_{(A,B,C,D)} := \gamma\mathcal{N}_{q^\dag}j^\dag_{(A,B,C,D)}$ and assume that noisy observations are given by
\begin{align}\label{20-6-17ct2}
\left( j^{(A,B,C,D)}_{\delta_{\ell}}, g^{(A,B,C,D)}_{\delta_{\ell}} \right) = \left( j^\dag_{(A,B,C,D)}+ \theta\cdot R_{j^\dag_{(A,B,C,D)}}, g^\dag_{(A,B,C,D)}+ \theta\cdot R_{g^\dag_{(A,B,C,D)}}\right) \quad \mbox{with} \quad \theta>0,
\end{align}
where $R_{j^\dag_{(A,B,C,D)}}$ and $R_{g^\dag_{(A,B,C,D)}}$ denote $\partial M^{h_\ell}\times 1$-matrices of random numbers on the interval $(-1,1)$.

With $\theta=0.1$ and $\ell=64$ the last line of Table \ref{b4} displays the numerical results for the case $(A,B,C,D)=(1,2,3,4)$ and $I=1$, which is repeated in the first line of Table \ref{b5} for comparison.

We now fix $D=4$. Let $(A,B,C)$ be equal to all permutations of the set $\{1,2,3\}$. Then, the equations \eqref{20-6-17ct1}--\eqref{20-6-17ct2} generate $I=6$ measurements. Similarly, let $(A,B,C,D)$ be all permutations of $\{1,2,3,4\}$ we get $I=16$ measurements. The numerical results for these two cases are presented in the two last lines of Table \ref{b5}, respectively.

Finally, in Figure \ref{h7} from left to right we show the computed numerical solution $q_\ell$ of the algorithm at the final iteration for $\ell=64$, $\theta=0.1$, i.e., $\delta_\ell=0.3308$, and $I=1,6,16$, respectively.

\begin{table}[H]
\begin{center}
\begin{tabular}{|c|l|l|l|l|l|}
\hline \multicolumn{6}{|c|}{ {\bf Numerical results for $\ell=64$, $\theta=0.1$} }\\
\hline
\mbox{Number of observations} $I$  &\scriptsize {\bf Iterate} &\scriptsize {\bf Tolerance} &\scriptsize $L^2_q$ &\scriptsize $L^2_\mathcal{N}$ &\scriptsize $L^2_\mathcal{D}$\\
\hline
1  & 1000& 2.2546e-2& 0.4059&  0.2067 &0.1077 \\
\hline
6 & 1000& 8.5684e-3 & 0.3159& 7.4901e-2& 4.4704e-2\\
\hline
16 & 1000 & 4.0133e-3& 0.2547& 5.6985e-2& 3.2211e-2\\
\hline
\end{tabular}
\caption{Numerical results for $\ell=64$, $\theta=0.1$, i.e., $\delta_\ell=0.3308$, and with multiple measurements $I=1,6,16$.}
\label{b5}
\end{center}
\end{table}

\begin{figure}[H]
\begin{center}
\includegraphics[scale=0.27]{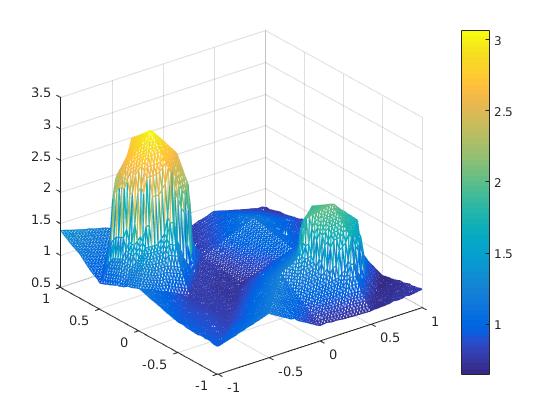}~~~~~~~~
\includegraphics[scale=0.27]{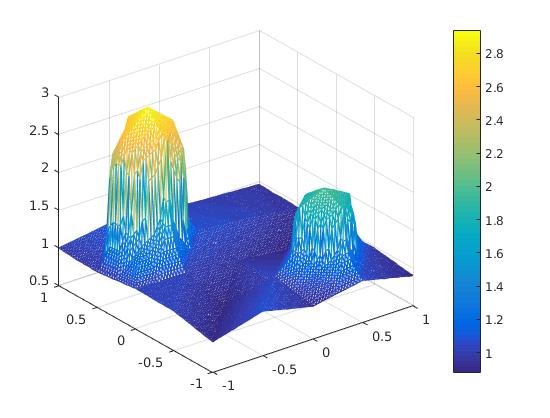} ~~~~~~~~
\includegraphics[scale=0.27]{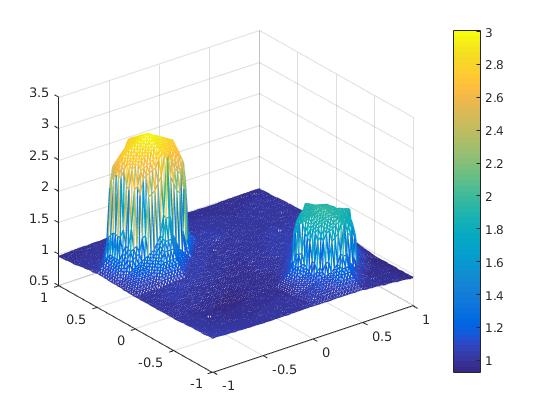}
\end{center}
\caption{Computed numerical solution $q_\ell$ of the algorithm at the final iteration for $\ell=64$, $\theta=0.1$, i.e., $\delta_\ell=0.3308$, and with multiple measurements $I=1,6,16$, respectively.}
\label{h7}
\end{figure}
We observe that the use of multiple measurements improves the solution to yield an acceptable result even in the presence of relatively large noise.
\end{example}

\section*{Acknowledgments} 

The authors M.~ Hinze, B.~ Kaltenbacher and T.N.T.~ Quyen would like to thank the referees and the editor for their valuable comments and suggestions which helped to improve our paper.

\end{document}